\documentclass[11pt]{amsart}

\usepackage{enumerate,graphics}
\usepackage[utf8]{inputenc}

\usepackage[%
  breaklinks=true, 
  colorlinks=true, 
  linkcolor=black, 
  anchorcolor=black,
  citecolor=black,
  filecolor=black,
  menucolor=black,
  pagecolor=black,
  urlcolor=blue,
  bookmarks,
  bookmarksnumbered,
  hyperfootnotes=false,
  pdfpagelabels,
  pdfstartview=FitH, 
  pdfpagemode=UseOutlines, 
  pdfnewwindow=true,
  pagebackref=false, 
  bookmarksopen=false,
  bookmarks=true
]{hyperref}
\usepackage[figure]{hypcap}

\newcommand{\hyp}{\mathbf{H}^3} 
\newcommand{\ball}{\mathbf{B}^3}
\newcommand{\Qthree}{\mathbf{Q}^3} 
\newcounter{primedtheoremdummy}
\newcommand{\remembercounter}[2]{\newcounter{#1}\setcounter{#1}{\value{#2}}}
\newcommand{\primedtheoremb}[1]{\setcounter{primedtheoremdummy}{\value{The}}\setcounter{The}{\value{#1}}\renewcommand{\theThe}{\arabic{The}'}}
\newcommand{\primedtheoreme}{\setcounter{The}{\value{primedtheoremdummy}}\renewcommand{\theThe}{\arabic{The}}}

\DeclareMathOperator{\End}{End}
\DeclareMathOperator{\Hom}{Hom}
\DeclareMathOperator{\GL}{GL}

\DeclareMathOperator{\SL}{SL}

\DeclareMathOperator{\SU}{SU}
\DeclareMathOperator{\Span}{Span}
\DeclareMathOperator{\Id}{Id}

\DeclareMathOperator{\gl}{gl}

\renewcommand{\Im}{\operatorname{Im}}
\renewcommand{\Re}{\operatorname{Re}}

\renewcommand{\ker}{\operatorname{ker}}
\newcommand{\im}{\operatorname{im}}
\DeclareMathOperator{\ord}{ord}

\newcommand{\R}{\mathbb{R}}
\newcommand{\C}{\mathbb{C}}
\newcommand{\N}{\mathbb{N}}

\renewcommand{\H}{\mathbb{H}}  

\newcommand{\CP}{\mathbb{CP}}
\newcommand{\HP}{\mathbb{HP}}

\renewcommand{\i}{\mathbf{i}}        

\renewcommand{\k}{\mathbf{k}}

\newcommand{\theset}[2]{\{\,#1\mid#2\,\}}
\newcommand{\dvector}[1]{{\left(\begin{matrix}#1\end{matrix}\right)}}
\newcommand{\tvector}[1]{{\left(\begin{smallmatrix}#1\end{smallmatrix}\right)}}

\newcommand{\tvectork}[1]{{\left[\begin{smallmatrix}#1\end{smallmatrix}\right]}}

\theoremstyle{plain}
\newtheorem{The}{Theorem}
\newtheorem*{The*}{Theorem}
\newtheorem{Pro}[The]{Proposition}
\newtheorem{Lem}[The]{Lemma}
\newtheorem{Cor}[The]{Corollary}

\theoremstyle{definition}
\newtheorem*{Def}{Definition}

\theoremstyle{remark}
\newtheorem{Rem}[The]{Remark}

\numberwithin{equation}{section}

\begin{document}

\title{Bryant Surfaces with Smooth Ends}

\author{Christoph Bohle}
\author{G.\ Paul Peters}

\address{Institute für Mathematik, 
Technische Universit{\"a}t Berlin,
Stra{\ss}e des 17.\ Juni 136,
10623 Berlin,
Germany}

\email{bohle@math.tu-berlin.de, peters@math.tu-berlin.de}


\date{\today}

\thanks{Both authors supported by DFG SPP 1154.}

\begin{abstract}
  A smooth end of a Bryant surface is a conformally immersed punctured disc of
  mean curvature $1$ in hyperbolic space that extends smoothly through the
  ideal boundary. The Bryant representation of a smooth end is well defined on
  the punctured disc and has a pole at the puncture. The Willmore energy of
  compact Bryant surfaces with smooth ends is quantized. It equals $4\pi$
  times the total pole order of its Bryant representation. The possible
  Willmore energies of Bryant spheres with smooth ends are $4\pi(\N^*\setminus
  \{2,3,5,7\})$. Bryant spheres with smooth ends are examples of soliton
  spheres, a class of rational conformal immersions of the sphere which also
  includes Willmore spheres in the conformal 3--sphere $S^3$. We give explicit
  examples of Bryant spheres with an arbitrary number of smooth ends. We
  conclude the paper by showing that Bryant's quartic differential $\mathcal Q$
  vanishes identically for a compact surface in $S^3$ if and only if it is the
  compactification of either a complete finite total curvature Euclidean
  minimal surface with planar ends or a compact Bryant surface with smooth
  ends.
\end{abstract}

\maketitle

\section{Introduction}

Surfaces of constant mean curvature $1$ in hyperbolic space have
attracted much attention since Bryant's fundamental paper \cite{Br87}.
They are nowadays called Bryant surfaces. Bryant surfaces are the
hyperbolic analogue to minimal surfaces in Euclidean space, because
their Gauss--Mainardi--Codazzi equations are virtually the same. This
leads to a local correspondence between Euclidean minimal and Bryant
surfaces which is known as the cousin relation.  Similar to Euclidean
minimal surfaces which are characterized by the holomorphicity of
their Gauss map, Bryant surfaces are characterized by the
holomorphicity of their hyperbolic Gauss map.  In contrast to the
Euclidean minimal case, the hyperbolic Gauss map of a complete finite total
curvature Bryant surface may not extend through the
ends. The ends at which the hyperbolic Gauss map extends
holomorphically are called regular ends, cf.~\cite{UY93,CHR01,Ro02}.

Smooth Bryant ends in hyperbolic space are a direct analog to planar
minimal ends in Euclidean space: both extend to immersions through the
ideal boundary of the space form, i.e., the 2--sphere or the point at
infinity, respectively. We show that a Bryant end is smooth if and
only if its Bryant representation $F$ is a holomorphic null immersion
of the punctured disc into $\SL(2,\C)$ with a pole at the puncture
such that $F'F^{-1}$ has a second order pole. In other words, a smooth
Bryant end is a regular end which is $\mathcal H^3$--reducible
\cite{RUY97,RUY04} (i.e., the representation $F$ has no monodromy)
such that $F'F^{-1}$ has a second order pole. The proof is based on a
theorem of Collin, Hauswirth, and Rosenberg \cite{CHR01} according to
which a properly embedded Bryant annular end is regular.  The
corresponding characterization in the Euclidean case is that a minimal
end is planar if and only if it is the real part of a holomorphic null
immersion into $\C^3$ with a simple pole at the end. In the case of
Bryant surfaces there are two types of smooth ends: those asymptotic
to horospheres and those asymptotic to smooth catenoid cousins in the
sense of \cite{ET01}. At a smooth horospherical end $F$ has a simple
pole and $F'F^{-1}$ automatically has a second order pole. At a smooth
catenoidal end $F$ has a pole of higher order.

Besides the cousin relation there is another correspondence between
Bryant surfaces and Euclidean minimal surfaces: if one interprets
$\SL(2,\C)$ as an affine representation of the complex 3--quadric
$Q^3\subset \CP^4$ and $\C^3$ as stereographic projection of $Q^3$
then both types of surfaces are represented by holomorphic null
immersions into $Q^3$.  The $Q^3$--valued representation of a smooth
Bryant end or a planar Euclidean minimal end extends to an
\emph{immersion} through the end. In case of smooth horospherical
Bryant ends and Euclidean minimal planar ends the holomorphic null
immersion transversally intersects the corresponding hyperplane at
infinity.  In case of smooth catenoidal Bryant ends the intersection
is non--transversal.

The analogy goes further if one considers the Willmore energy of
surfaces in the conformal 3--sphere $S^3=\R^3\cup \{\infty\}$ (see the
appendix for the relation of the Willmore energy and the total
absolute curvature). As proven by Bryant~\cite{Br84}, all complete
finite total curvature minimal surfaces with planar ends in $\R^3$
extend to compact Willmore surfaces in $S^3$, i.e., critical points of
the Willmore energy, and all Willmore \emph{spheres} in $S^3$ are
obtained this way. Bryant surfaces with smooth ends are not Willmore,
but constrained Willmore~\cite{BPP1}. Nevertheless, they obey the same
quantization of the Willmore energy as Euclidean minimal surfaces with
planar ends: in both cases the Willmore energy is $4\pi d$ where
$d\in\N$ is the total pole order of the corresponding null curve in
$\SL(2,\C)$ or $\C^3$, respectively. As in the case of Euclidean
minimal spheres with planar ends~\cite{Br84,Br88}, the
possible Willmore energies of Bryant \emph{spheres} with smooth ends
are $4\pi(\N^*\setminus\{2,3,5,7\})$. In particular, we show that 
Bryant spheres with exactly $d$ horospherical smooth ends exist if and
only if $d\in\N^*\setminus\{2,3,5,7\}$. There is no
restriction on the number of ends if one allows for catenoidal
ends: we give explicit rational conformal parametrizations for spheres
with an arbitrary number of smooth ends.

\begin{figure}[t]
  \centering
  \scalebox{.6}{\includegraphics{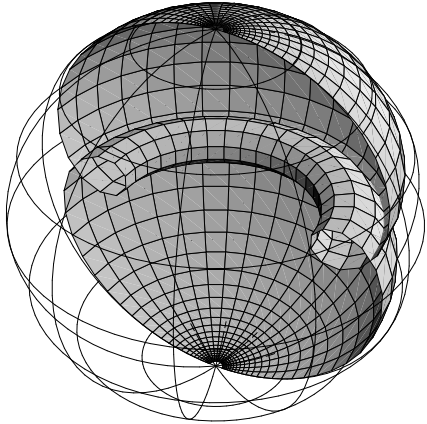}}
  \caption{Dirac sphere alias catenoid cousin.}
  \label{fig:Dirac1}
\end{figure}

We were led to consider Bryant surfaces with smooth ends during our
investigations of soliton spheres \cite{Pe04,BoPe} when we observed
that the simplest non--trivial Dirac sphere \cite{Ri97}, a surface of
revolution related to a 1--soliton solution of the mKdV
equation~\cite{Ta99}, is a catenoid cousin with smooth ends, see
Figure~\ref{fig:Dirac1}.  We prove that all Bryant spheres with smooth
ends are soliton spheres. The analogous result that all Euclidean
minimal spheres with planar ends are soliton spheres is proven
in~\cite{Pe04,BoPe}.

In the last section we give a uniform Möbius geometric characterization of
complete Euclidean minimal surfaces with finite total curvature and planar
ends and compact Bryant surfaces with smooth ends in terms of Bryant's
quartic differential~$\mathcal Q$: we show that $\mathcal Q$ vanishes
identically for an immersion of a compact surface into $S^3$ if and only if it
is the compactification of either a Euclidean minimal surface with planar ends
or a compact Bryant surface with smooth ends.

\subsection*{Acknowledgments}

We thank Ulrich Pinkall, Alexander Bobenko, Tatya\-na Pavlyukevich, Udo
Hertrich--Jeromin, and Wayne Rossmann for helpful discussions.


\section{The Bryant Representation}\label{sec:bryant_representation}

The unit 3--ball $\ball=\theset{x\in\R^3}{|x|<1}$ with the metric
$ds^2=\frac{4|dx|^2}{(1-|x|^2)^2}$ is the Poincar{\'e} ball model of hyperbolic
3--space. A surface in $\ball$ is called a \emph{Bryant surface} if it has
constant mean curvature one. A fundamental property of Bryant surfaces is that
they posses a so called \emph{Bryant representation}~\cite{Br87} in terms of
holomorphic data similar to the Weierstrass representation of minimal surfaces
in $\R^3$: a conformal immersion $f\colon M\to\ball$ of a Riemann surface $M$
parametrizes a Bryant surface if and only if there exists a holomorphic
null immersion $F\colon \tilde M\to\SL(2,\C)$ (null meaning $\det F' =0$) defined on
the universal covering $\tilde M$ of $M$ such that
\begin{equation}
  \label{eq:Formula_in_Ball}
  f=\frac{(x_1,x_2,x_3)}{x_0+1},\qquad \textrm{where} \qquad
 \dvector{x_0+x_3&x_1+x_2{\i}\\x_1-x_2{\i}&x_0-x_3}
 :=F\bar F^t.
\end{equation}
The holomorphic null immersion $F$ representing the Bryant surface $f$ is
unique up to right multiplication by a constant $\SU(2)$ matrix. Left
multiplication of $F$ with $\SL(2,\C)$ matrices yields all congruent
Bryant surfaces.

Although the ball model of hyperbolic space is best suited for the definition
of smooth ends, the investigation of the ends turns out to be much simpler in
the half space model $\hyp=\theset{(x_1,x_2,x_3)\in\R^3}{x_3>0}$ with the
metric $ds^2=\frac{|dx|^2}{x_3^2}$. In the half space model the Bryant surface
corresponding to a holomorphic null immersion $F= \tvector{a&b\\c&d}\colon
\tilde M\to \SL(2,\C)$ is given by 
\begin{align}
  \label{eq:Formula_in_halfspace}
  x_1+{\i} x_2&=\frac{a\bar c+b\bar d}{|c|^2+|d|^2},&
  x_3&=\frac1{|c|^2+|d|^2}.
\end{align}
These formulas may be derived from \eqref{eq:Formula_in_Ball} by
applying the orientation preserving isometry: $\ball\to\hyp$,
$(x_1,x_2,x_3)\mapsto2\frac{(x_1,x_2,1-x_3)}{|(x_1,x_2,1-x_3)|^2}-(0,0,1)$.

$F$ is a holomorphic null immersion into $\SL(2,\C)$ if and only if
$F^{-1}$ is a holomorphic null immersion. The surface corresponding to
$F^{-1}$ is called the \emph{dual} Bryant surface of the one
represented by $F$.

\section{Smooth Ends}\label{sec:smooth_ends}

In this section we characterize Bryant representations of smooth
Bryant ends. We reformulate this characterization in terms of
holomorphic null immersions into $Q^3$.  This allows to determine the
possible Willmore energies of Bryant spheres with smooth ends.

\begin{Def}
  We call a Bryant surface $E$ in the Poincar{\'e} ball model
  $\ball\subset\R^3$ of hyperbolic space a \emph{smooth Bryant end} if
  there is a point $p_\infty\in \partial \ball$ on the asymptotic
  boundary such that $E\cup\{p_\infty\}$ is a conformally immersed
  open disc in $\R^3$.  We call a Bryant surface a \emph{compact
    Bryant surface with smooth ends} if it is conformally equivalent
  to a compact Riemann surface minus a finite number of points such
  that each of these points has a punctured neighborhood that is a
  smooth Bryant end. We call it a \emph{Bryant sphere with
    smooth ends} if the underlying compact Riemann surface is the
  sphere.
\end{Def}

The definition is invariant under hyperbolic isometries and one might
replace the ball model $\ball$ by the half space model $\hyp$ as long
as no end goes to $\infty$, because hyperbolic isometries extend to
M{\"o}bius transformations of the conformal 3--sphere
$S^3=\R^3\cup\{\infty\}$.

We prove now that the Bryant representation of a smooth end is well
defined on the punctured unit disc $\Delta^*=\Delta\setminus\{0\}$ and
has a pole at zero.  The proof is based on a result of Collin,
Hauswirth, and Rosenberg~\cite{CHR01} about the behavior of the
holomorphic null immersion $F$ corresponding to a properly embedded
Bryant annular end.

\remembercounter{Tsmoothend}{The}
\begin{The}\label{T:smooth_end}
  If $E$ is a smooth Bryant end, then the corresponding null immersion
  $F\colon \tilde \Delta^*\to\SL(2,\C)$ is well defined on $\Delta^*$, has a
  pole at zero, and $F'F^{-1}$ has a pole of order 2 at zero. Conversely, if
  $F\colon \Delta^*\to\SL(2,\C)$ is a holomorphic null immersion with a pole
  at zero such that $F'F^{-1}$ has a pole of order 2 at zero, then the Bryant
  surface corresponding to $F$ is a smooth Bryant end.
\end{The}

\begin{Rem}
  More precisely, Lemma~\ref{L:poles_of_omega} implies that if
  the holomorphic null immersion $F$ has a pole of order $1$ at zero,
  then the Bryant end is asymptotic, in the sense of \cite{ET01}, to
  the end of a horosphere, and $F'F^{-1}$ has automatically a pole of
  order $2$ at zero. If $F$ has a pole of order $n>1$, then the end is
  asymptotic to the end of a smooth catenoid cousin with Bryant
  representation
  \begin{equation}
    \label{eq:CC_Formula}
    F=\frac1{\sqrt{2\mu+1}}\tvector{(\mu+1)z^\mu &\mu z^{-(\mu+1)}\\
      \mu z^{\mu+1}&(\mu+1)z^{-\mu}},
  \end{equation}
  where $\mu=n-1$, cf.~\cite{Br87}. The fact that the catenoid cousin
  represented by~(\ref{eq:CC_Formula}) has smooth ends if and only if
  $\mu=n-1$ can be easily derived
  from~(\ref{eq:Formula_in_halfspace}). A similar observation for
  trinoids is mentioned in \cite{BPS02}.
\end{Rem}

\begin{proof}
  Let $E\subset\ball$ be a smooth Bryant end, $f\colon \Delta\to\R^3$
  a conformal immersion such that $f_{|_{\Delta^*}}$ parametrizes $E$,
  and $F\colon \tilde \Delta^*\to\SL(2,\C)$ the corresponding
  holomorphic null immersion. Then there is a closed disc $\bar
  \Delta_r\subset\Delta$ such that $f_{|_{{\bar \Delta_r}^*}}$ is a
  proper embedding into $\ball$.  In \cite{CHR01} it is shown in the
  proofs of Theorem~3 and 4 that $F$ on $\tilde \Delta_r^*$ is, up to
  right--multiplication by a constant $\SU(2)$ matrix, the product of
  a holomorphic $\SL(2,\C)$--valued map defined on $\Delta^*_r$ with a
  pole at $z=0$ and $\tvector{z^{-\nu}&0\\0&z^\nu}$ for some
  $\nu\in\R$. (One may also derive this, using the description of the
  behavior of $F$ at regular ends given in \cite{UY93}, from the fact
  that properly embedded Bryant annular ends are regular
  \cite{CHR01}.) Hence there is a holomorphic map $\tilde F$ on
  $\Delta_r$ that does not vanish at zero, $\alpha\in
  [-\frac12,\frac12]$, and $n\in\N$ such that
  \begin{equation*}
    F(z)=z^{-n}\tilde F(z)\dvector{z^{-\alpha }&0\\0&z^\alpha},
  \end{equation*}
  for all $z\in\Delta^*_r$.
  
  We now show that smoothness of $f$ at zero implies
  $\alpha=0$. Suppose that $\alpha\neq 0$ and let $\tilde
  F=\tvector{a&b\\c&d}$.  Multiplying $F$ by a constant $\SL(2,\C)$
  matrix from the left and a constant $\SU(2)$ matrix from the right,
  we may assume that $c(0)\neq0$ and $a(0)=0$. The conformal immersion
  into $\hyp$ satisfies
  \begin{align*}
    x_1+{\i} x_2&=\frac{a\bar c+|z|^{4\alpha}b\bar d}{|c|^2+|z|^{4\alpha}|d|^2},&
    x_3&=\frac{|z|^{2n+2\alpha}}{|c|^2+|z|^{4\alpha}|d|^2}.
  \end{align*}
  
  Thus $x_1$, $x_2$, and $x_3$ tend to zero if $z$ tends to zero.  Therefore,
  for $f$ to be smooth in particular $x_3$ has to be smooth at zero. If
  $0<\alpha<\frac12$, then $x_3(z)=o( |z|^{2n})$ but not
  $O(|z|^{2n+1})$ as $z\to0$ and by Taylor's theorem $x_3$ could not be
  $C^{2n+1}$. Similarly, if $-\frac12 < \alpha<0$, then $x_3(z)=o(
  |z|^{2n-1})$ but not $O(|z|^{2n})$ as $z\to0$ and $x_3$ could not be
  $C^{2n}$.  For $|\alpha|=\frac12$ the denominator of $x_3$ would be smooth
  and therefore, because the denominator does not vanish at zero and the
  numerator $|z|^{2n\pm1}$ is not smooth, $x_3$ could not be smooth.  Thus
  $\alpha$ has to be zero.

  The claim about the pole order of $F'F^{-1}$ and the converse follows from
  Lemma~\ref{L:poles_of_omega} below.
\end{proof}

We now derive a normal form for the Bryant representation $F$ at a
pole.

\begin{Lem}[Normal form]\label{L:normal_form}
  Let $F\colon \Delta^*\to\SL(2,\C)$ be a holomorphic map with a pole of order
  $n\in\N^*$ at zero. Then there exist matrices $A\in\SL(2,\C)$, $B\in
  \SU(2)$, and holomorphic functions $a,b,c,d\colon \Delta\to\C$ such that
  \begin{align*}
    AFB&=z^{-n}\dvector{a&b\\c&d},& a(0)&=a'(0)=b(0)=c(0)=0,& d(0)&\neq0.
  \end{align*}
  If $F$ is null, then the vanishing orders of  $a$ and $bc$ at
  $z=0$ satisfy
  \begin{align*}
    \ord_0(a)&\geq2n,& \ord_0(bc)&=2n.
  \end{align*}
\end{Lem}

In the half plane model \eqref{eq:Formula_in_halfspace}, the
fact that $F$ is in normal form implies that the corresponding Bryant end
converges to $0\in \partial\hyp$.

\begin{proof}
  Since $F$ has a pole of order $n$ at zero there exist holomorphic maps
  $a,b,c,d\colon\Delta\to\C$ such that
  \begin{equation*}
    F=z^{-n}\dvector{a&b\\c&d},
  \end{equation*}
  and one of the functions $a,b,c,d$ does not vanish at zero. Multiplying $F$
  by $\tvector{0&-1\\1&0}$ from the left or right we may assume that
  $d(0)\neq0$.  Multiplying $F$ from the right by
  \begin{equation*}
    \left(1+\tfrac{|c(0)|^2}{|d(0)|^2}\right)^{-\frac12}
    \dvector{1&\frac{\bar c(0)}{\bar d(0)}\\-\frac{c(0)}{d(0)}&1}
  \end{equation*}
  one gets $c(0)=0$ while $d(0)\neq0$ is preserved. Multiplying $F$ from the
  left by $\tvector{1&-\frac{b(0)}{d(0)}\\0&1}$ one obtains $b(0)=0$ while
  $c(0)=0$ and $d(0)\neq0$ is preserved. Since $F$ has determinant one (and we
  did not change this by our multiplications), we have $ad-bc=z^{2n}$.  Thus
  $a$ vanishes to the second order at zero, because $n\geq1$,
  $d(0)\neq0$, and $bc$ vanishes to
  the second order at zero.
  
  If $F$ is null, then $\det F'=0$ and $ad-bc=z^{2n}$ imply
  $a'd'-b'c'=n^2z^{2n-2}$. Suppose $\ord_0(bc)=k<2n$, then $\ord_0(a)=k$. Hence
  $\ord_0(a'd')\geq k-1$ and $\ord_0(b'c')=k-2<2n-2$, which is a contradiction
  to $a'd'-b'c'=n^2z^{2n-2}$. If $\ord_0(bc)=k>2n$, then $\ord_0(a)=2n$.
  Hence $\ord_0(a'd')\geq2n-1$ and $\ord_0(b'c')=k-2>2n-2$ which again
  contradicts $a'd'-b'c'=n^2z^{2n-2}$. Hence $\ord_0(bc)=2n$ and
  $\ord_0(a)\geq2n$.
\end{proof}

\begin{Lem}\label{L:poles_of_omega}
  Let $F\colon\Delta^*\to\SL(2,\C)$ be a holomorphic null immersion with a
  pole at zero and let $E$ be the Bryant surface corresponding to
  $F$.
  \begin{enumerate}[(i)]
  \item $E$ is a smooth end if and only if $F'F^{-1}$ has a pole of order $2$
    at zero.
  \item The dual Bryant surface of $E$ is a smooth end if and only if
    $F^{-1}F'$ has a pole of order $2$ at zero.
  \end{enumerate}
  Both $E$ and its dual are smooth ends if and only if $F$ has a pole of order
  1.
\end{Lem}

\begin{proof}
  Let $n\in\N^*$ be the pole order of $F$ and assume that $F$ has the normal
  form of Lemma~\ref{L:normal_form}. By \eqref{eq:Formula_in_halfspace}, $E$
  has the conformal parametrization
  \begin{align*}
    x_1+{\i} x_2&=\frac{a\bar c+b\bar d}{|c|^2+|d|^2},&
    x_3&=\frac{|z|^{2n}}{|c|^2+|d|^2}
  \end{align*} 
  defined on the punctured disc $\Delta^*$. The parametrization can be
  smoothly extended to zero by setting
  $(x_1(0),x_2(0),x_3(0))=0\in\partial\hyp$.
  
  By Lemma~\ref{L:normal_form}, $a$ vanishes to the second order at zero and
  we have
  \begin{align*}
    dx_1+{\i} d x_2 &= \frac{b'\bar d}{|c|^2+|d|^2} dz+ O( |z|),&&\text{as $z\to
      0$}.
  \end{align*}
  Hence $b'(0)\neq0$ if and only if $(x_1,x_2,x_3)\colon \Delta \to \R^3$ is
  an immersion, i.e., if and only if $E$ is a smooth end.  Because, up to
  isometry, passing to the dual surface amounts to interchanging $c$ and $b$,
  $c'(0)\neq0$ if and only if the dual surface is a smooth end.  By
  Lemma~\ref{L:normal_form} we have $\ord_0(bc)=2n$, therefore the case that
  both $b'(0)\neq 0$ and $c'(0)\neq0$ is equivalent to $n=1$.
  
  To complete the proof we have to show that $b'(0)\neq 0$ is equivalent to
  $F'F^{-1}$ having a pole of order $2$ and $c'(0)\neq0$ is equivalent to
  $F^{-1}F'$ having a pole of order $2$. We have
  \begin{align*}
    F'F^{-1}=-nz^{-1}\Id  + z^{-2n}\dvector{a'd-b'c&-a'b+b'a\\c'd-d'c&-c'b+d'a}.
  \end{align*}
  Lemma~\ref{L:normal_form} implies $\ord_0(a'd-b'c)\geq 2n-1$,
  $\ord_0(-a'b+b'a)\geq 2n$, $\ord_0(c'd-d'c)= \ord_0(c)-1$, and
  $\ord_0(-c'b+d'a)\geq 2n-1$. Hence $F'F^{-1}$ has a pole of order $2$ if and
  only if $\ord_0(c)=2n-1$, which is, by Lemma~\ref{L:normal_form}, equivalent
  to $b'(0)\neq0$. Because passing to the dual surface essentially
  interchanges both the roles of $b$ and $c$ and the roles of $F'F^{-1}$ and
  $F^{-1} F'$, this also proves the statement for $c'(0)\neq0$.
\end{proof}

The following corollary to Theorem~\ref{T:smooth_end} characterizes the Bryant
representations of compact Bryant surfaces with smooth ends.

\remembercounter{Cbryantsurfacewithsmoothends}{The}
\begin{Cor}\label{C:bryant_surface_with_smooth_ends}
  The Bryant representation $F$ of a compact Bryant surface with
  smooth ends is a meromorphic null immersion with $\SU(2)$--monodromy
  of the underlying compact Riemann surface $M$ into $\SL(2,\C)$. The
  poles of $F$ correspond to the ends of the surface.  In particular,
  a Bryant sphere with smooth ends is represented by a rational
  $\SL(2,\C)$--valued null immersion.  Conversely, if $F$ is an
  $\SL(2,\C)$--valued null immersion with $SU(2)$--mo\-no\-dro\-my of
  a compact Riemann surface $M$ such that $F'F^{-1}$ has only poles of
  order 2, then $F$ represents a compact Bryant surface with smooth
  ends.
\end{Cor}

\begin{figure}[h]
  \centering 
  \scalebox{.6}{\includegraphics{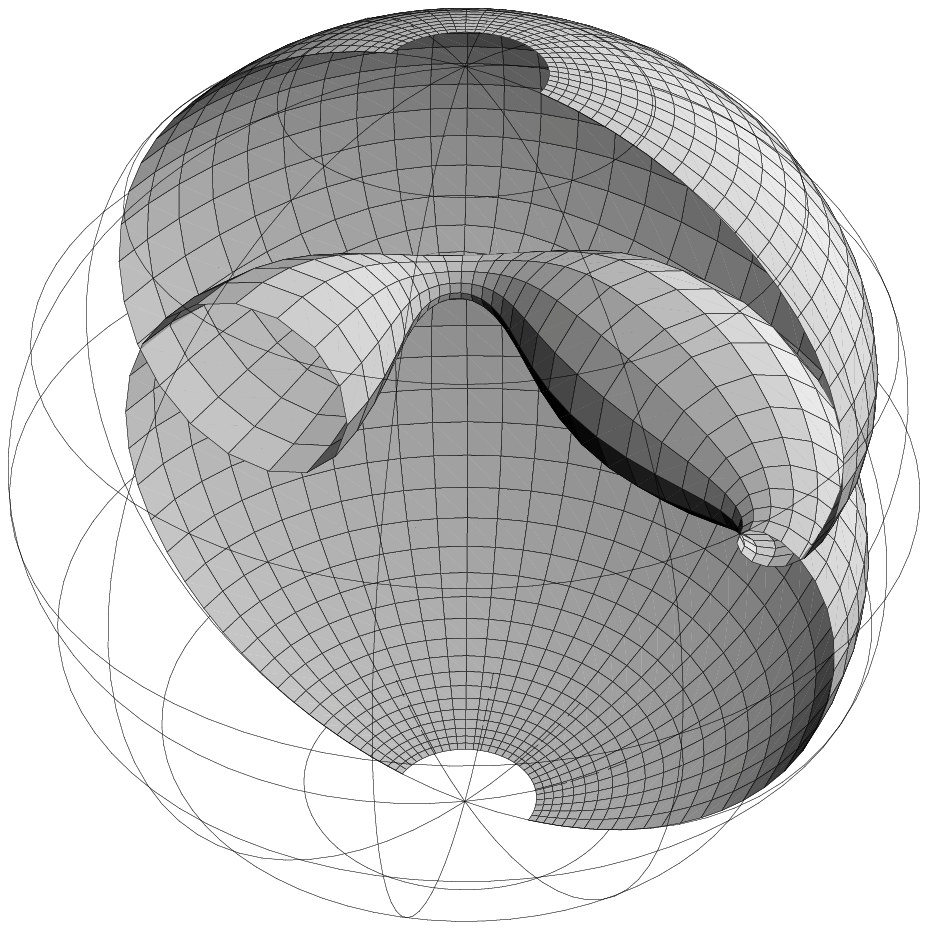}}\hfill
  \scalebox{.6}{\includegraphics{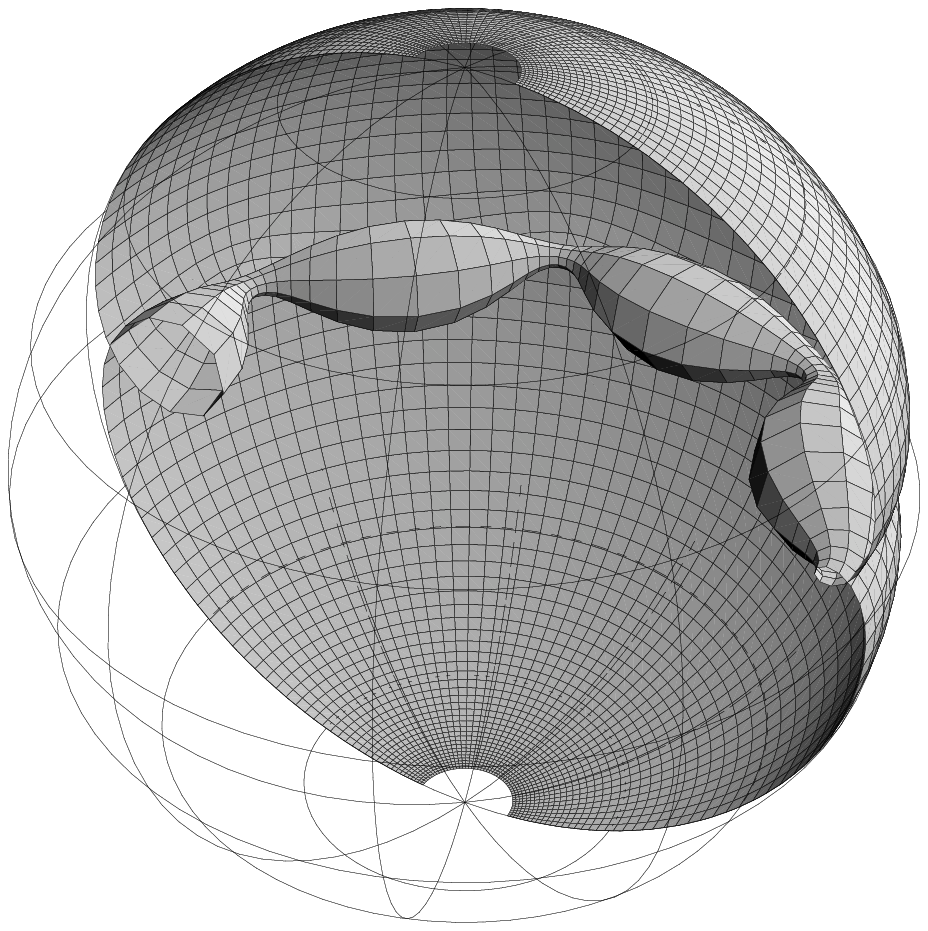}}\hfill
  \caption{Warped catenoid cousins with smooth ends.}
  \label{fig:warped_cc}
\end{figure}

The corollary implies that Bryant spheres with smooth ends are
complete $\mathcal H^3$--reducible Bryant spheres \cite{RUY97,RUY04}
of finite total curvature with regular ends, and that all its
$\mathcal H^3$--deformations are Bryant spheres with smooth
ends. Examples of the $\mathcal H^3$--deformation are the warped
catenoid cousins \cite{RUY03} represented by $FA$, where $F$ is the
representation of a catenoid cousin with smooth ends as
in~\eqref{eq:CC_Formula} and $A\in \SL(2,\C)$.
Figure~\ref{fig:warped_cc} shows two spheres obtained for
$A=\tvector{1&1/2\\0&1}$ and $\mu=1$ or $3$.

A more geometric approach  to the study of meromorphic null immersions
into $\SL(2,\C)$ is to interpret $\SL(2,\C)$ as the affine part
$\Qthree\cap\{e\neq0\}$ of the quadric
\begin{equation*}
  \Qthree=\theset{[a,b,c,d,e]\in\CP^4}{ad-bc-e^2=0}\subset\CP^4.
\end{equation*}
More precisely,
\begin{align*}
  \SL(2,\C)&\overset{1:1}{\longleftrightarrow}
  \Qthree\cap\theset{[a,b,c,d,e]\in\CP^4}{e\neq0}\\ 
  F=\tvector{a&b\\c&d} &\longleftrightarrow \Phi=[a,b,c,d,1].
\end{align*}
A holomorphic map $\Phi=[a,b,c,d,e]$ into $\Qthree$ is null, i.e.,
$a'd'-b'c'-e'^2=0$, if and only if the corresponding
$\SL(2,\C)$--valued meromorphic map $F=\frac1e\tvector{a&b\\c&d}$ is
null. By \eqref{eq:Formula_in_halfspace}, the Bryant surfaces that
corresponds to a holomorphic null immersion $\Phi=[a,b,c,d,e]\colon
M\to \Qthree$ is
\begin{align*}
  x_1+{\i} x_2&=\frac{a\bar c+b\bar d}{|c|^2+|d|^2},&
  x_3&=\frac{|e|^2}{|c|^2+|d|^2}.
\end{align*}
The zeroes of $e$ are the ends of the Bryant surface and the pole
order of $F$ at the ends coincides with the intersection order of
$\Phi$ with the hyperplane $\{e=0\}$ at the end. This interpretation
of meromorphic null maps into $\SL(2,\C)$ was suggested to the authors
by Ulrich Pinkall.

The following proposition allows to  restate Theorem~\ref{T:smooth_end} and
Corollary~\ref{C:bryant_surface_with_smooth_ends} in terms of $\Phi$.

\begin{Pro}\label{P:Phi_immersed} 
  A holomorphic null map $\Phi\colon \Delta\to \Qthree$ for which
  $\Phi(0)$ is contained in the hyperplane $\{e=0\}$ is immersed at
  $0$ if and only if $F'F^{-1}$ or $F^{-1}F'$ has a pole of order $2$
  at $0$. The intersection of $\Phi$ with $\{e=0\}$ is transversal if
  and only if both $F'F^{-1}$ and $F^{-1}F'$ have poles of order $2$
  at $0$.
\end{Pro}

\begin{proof}
  If the holomorphic map $\Phi\colon \Delta\to \Qthree$ satisfies
  $\Phi(0)\in \{e=0\}$ then the corresponding $F\colon \Delta^*\to
  \SL(2,\C)$ has a pole at zero and there are holomorphic functions
  $a,b,c,d\colon \Delta\to\C$ one of which does not vanish at zero and
  $n\in\N^*$ such that $F=z^{-n}\tvector{a&b\\c&d}$ and
  $\Phi=[a,b,c,d,z^{n}]$. We may assume that $F$ is in normal form of
  Lemma~\ref{L:normal_form}, because replacing $F$ by $AFB$ does not
  change the pole orders of $F'F^{-1}$ or $F^{-1}F'$ and $\Phi$
  changes by a projective transformation that preserves
  $\{e=0\}$. Then, $a$ vanishes at least to the second order and $d$
  does not vanish at zero, so $\Phi$ is immersed if and only if
  $n=1$, $b'(0)\neq0$, or $c'(0)\neq0$. As we have seen in the proof of
  Lemma~\ref{L:poles_of_omega}, this is equivalent to $F'F^{-1}$
  or $F^{-1}F'$ having a pole of order $2$ at zero.
\end{proof}

Theorem~\ref{T:smooth_end} and
Corollary~\ref{C:bryant_surface_with_smooth_ends} can now be restated as
follows: \primedtheoremb{Tsmoothend}
\begin{The}\label{T:smooth_end_prime}
  If $E$ is a smooth Bryant end, then the corresponding null immersion
  $\Phi\colon \tilde \Delta^*\to\Qthree$ is well defined on $\Delta^*$ and
  extends to a holomorphic immersion on $\Delta$. Conversely, if $\Phi\colon
  \Delta\to\Qthree$ is a holomorphic null immersion such that $\Phi(0)$ is
  contained in the hyperplane $\{e=0\}$, then the corresponding Bryant surface
  or its dual Bryant surface is a smooth Bryant end.
\end{The}
\primedtheoreme

The generic case, i.e., $\Phi$ intersects the hyperplane $\{e=0\}$
transversally, is equivalent to both surfaces being smooth
horospherical Bryant ends.

\primedtheoremb{Cbryantsurfacewithsmoothends}
\begin{Cor}\label{C:bryant_surface_with_smooth_ends_prime}
  A compact Bryant surface with smooth ends is represented by a holomorphic
  null immersion $\Phi$ of the universal covering $\tilde M$ of a compact
  Riemann surface $M$ into $\Qthree$. In particular, a Bryant sphere with
  smooth ends is represented by a rational null immersion into $\Qthree$.
\end{Cor}
\primedtheoreme

The fact that Bryant spheres with smooth ends are represented by
rational null immersions into $\Qthree$ provides a strong link to
Willmore spheres in $S^3$, because these as well are related to
rational null immersion into $\Qthree$, see \cite{Br84}.  As in the
case of Willmore spheres, the Willmore energy of a Bryant sphere with
smooth ends is given by $4\pi \deg\Phi$. This follows from
Theorem~\ref{T:total_curvature} below, because $\deg\Phi$ is the total
pole order of $F$. Using Bryant's result~\cite{Br88} that the possible
degrees $d$ of rational null immersions into $\Qthree$ are the numbers
$d\in\N^*\setminus\{2,3,5,7\}$, we obtain the following theorem.

\begin{The}\label{T:gaps_for_spheres}
  The possible Willmore energies of Bryant spheres with smooth ends are $4\pi
  d$ with $d\in\N^*\setminus\{2,3,5,7\}$.
\end{The}

\begin{Cor}
  Bryant spheres with $d$ smooth ends that are all horospherical
  exist if and only if $d\in\N^*\setminus\{2,3,5,7\}$.
\end{Cor}

\begin{figure}[h]

  \centering
  \scalebox{.6}{\includegraphics{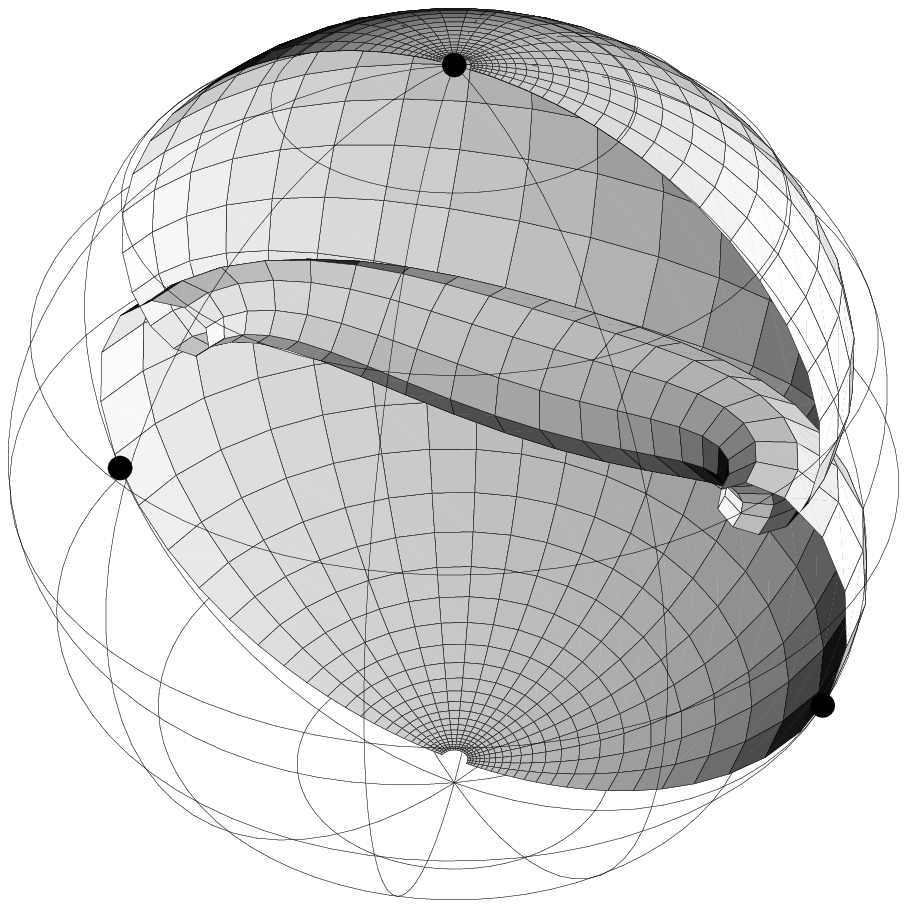}}\hspace{2cm}
  \scalebox{.6}{\includegraphics{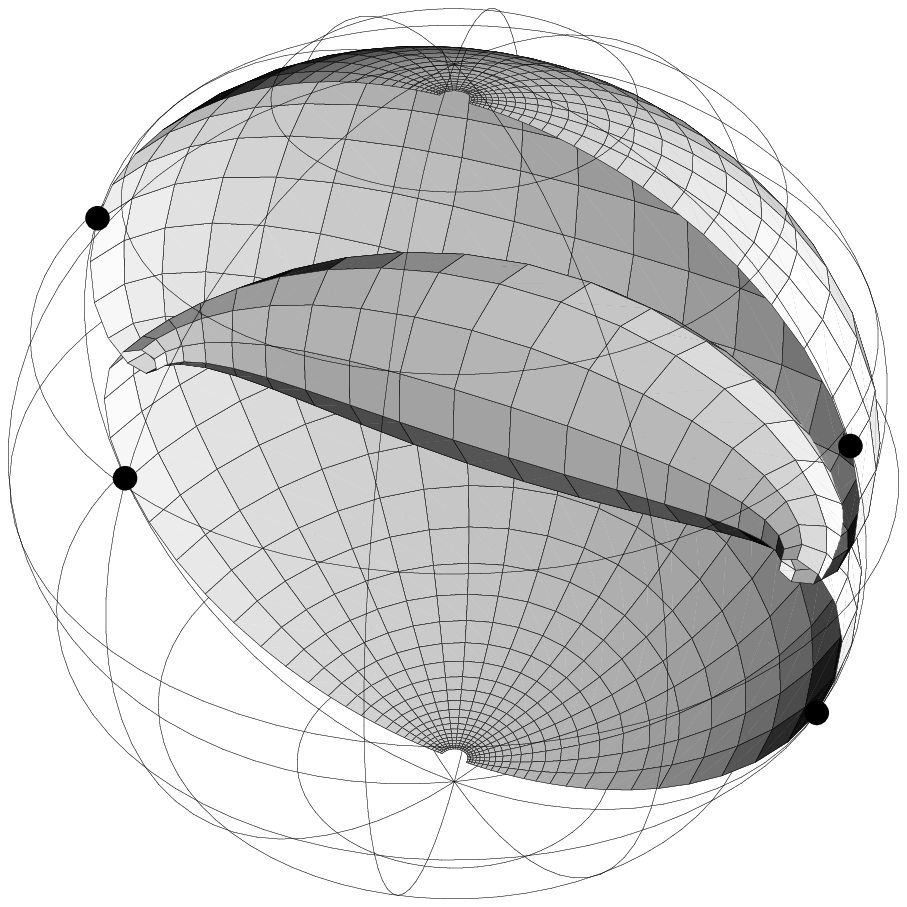}}\vspace{6ex}
  \caption{$3$-- and $4$--noid with smooth ends
  (marked points) and Willmore energy $16\pi$.}
  \label{fig:cc-deformation1}
\end{figure}

The generic case of horospherical smooth ends does not occur among the
rotationally symmetric examples of the catenoid cousins. However, generic
projective transformations of $\CP^4$ that preserve $\Qthree$ deform
catenoidal smooth ends into horospherical smooth ends.  As an example, we
apply the projective transformation $a\mapsto a$, $d\mapsto d$,
$\tvector{e&b\\c&-e}\mapsto A^{-1}\tvector{e&b\\c&-e}A$,
$A=\tvector{1&0\\-t&1}\tvector{1&s\\0&1}$ to the holomorphic null immersion
$\Phi=[a,b,c,d,e]$ of a catenoid cousin.  If $\mu=n-1$, $n\in
\N^*\setminus\{1\}$ in~\eqref{eq:CC_Formula} and
$st\in\C^*\setminus\{1,\frac1{2\mu+2},\frac{2\mu+1}{2\mu+2}\}$ one obtains a
Bryant sphere with $2n$ horospherical smooth ends and Willmore
energy~$8n\pi$. The surface on the right in Figure~\ref{fig:cc-deformation1}
is an example corresponding to $n=2$, $s=t=0.2$.

\begin{figure}[h]
  \centering
  \scalebox{.6}{\includegraphics{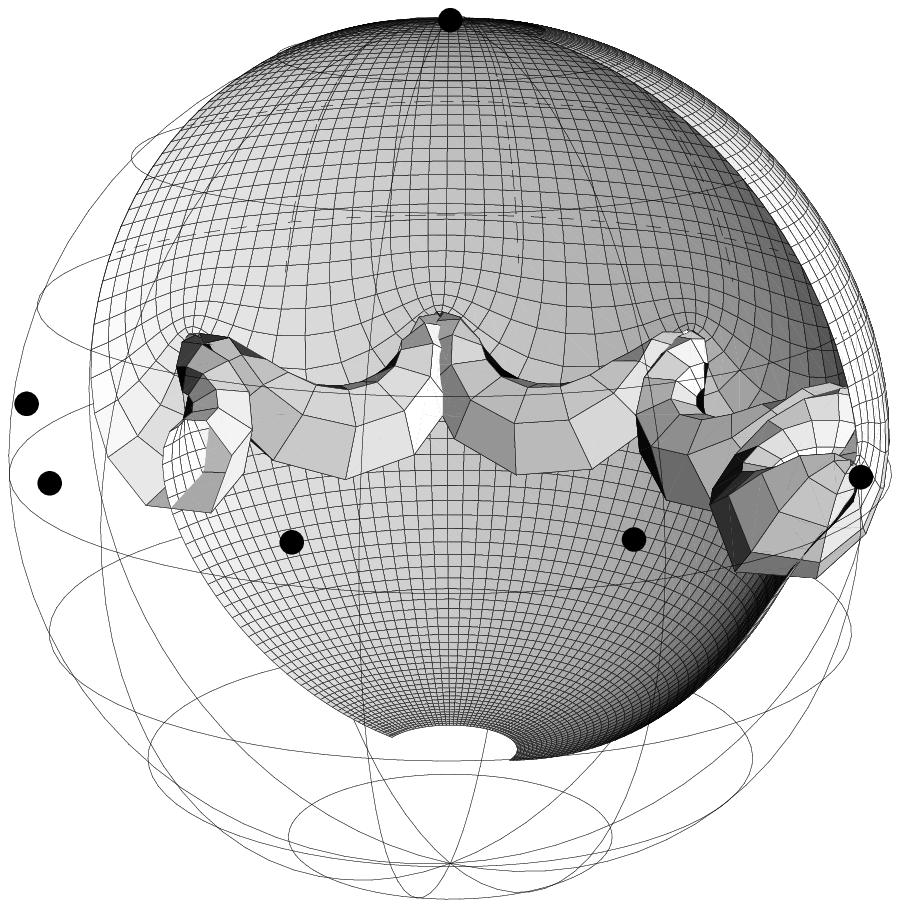}}\hspace{2cm}
  \scalebox{.6}{\includegraphics{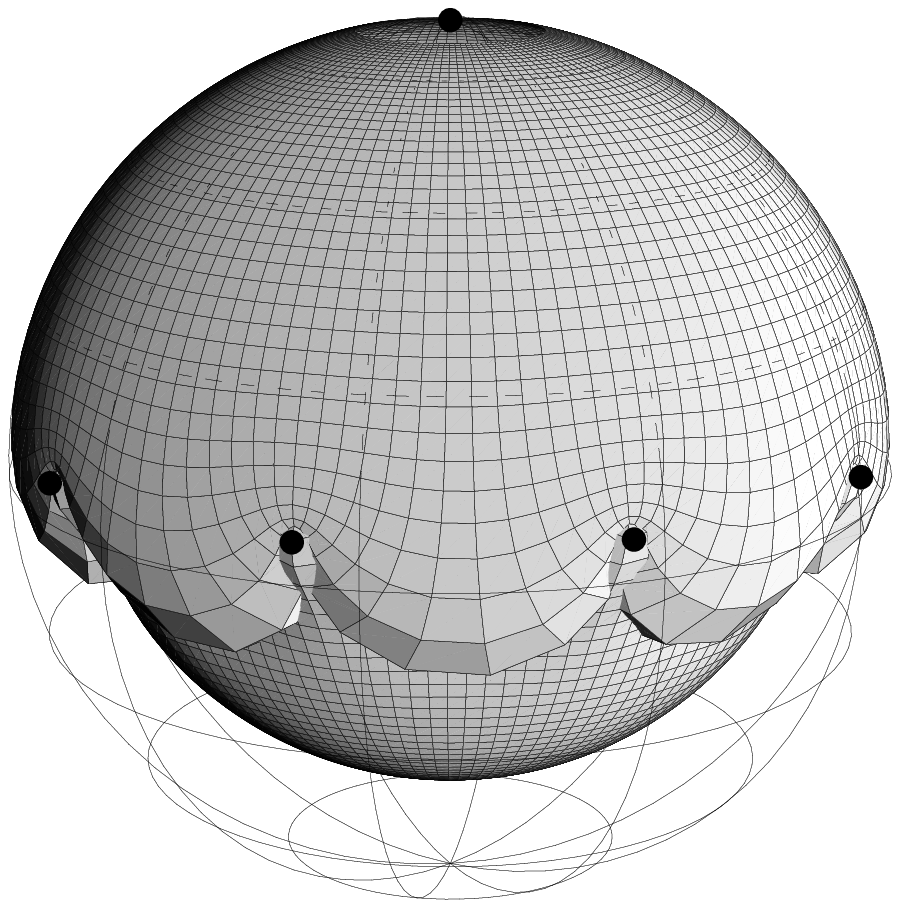}}\vspace{6ex}
  \caption{$10$--noid with smooth ends
  (marked points) and Willmore energy $72\pi$.}
  \label{fig:cc-deformation2}
\end{figure}

For $s\in \C^*$ and $t=0$ one gets Bryant spheres with 1
catenoidal and $n$ horospherical ends and Willmore energy $8n\pi$. The
left surface in Figure~\ref{fig:cc-deformation1} has 3 smooth ends and
corresponds to the parameters $n=2$, $s=0.2$, and $t=0$, and
Figure~\ref{fig:cc-deformation2} shows two views of the Bryant spheres
with 10 smooth ends corresponding to $n=9$, $s=.4$, and $t=0$.

The explicit rational conformal immersions obtained by deforming catenoid
cousins with smooth ends show that \emph{there exist Bryant spheres with an
  arbitrary number of smooth ends.}


\section{Bryant Surfaces as Darboux Transforms \\ of the Round Sphere --- a Quaternionic Approach}\label{sec:darboux}

This section is devoted to a M{\"o}bius geometric interpretation of
the Bryant representation due to Hertrich-Jeromin, Musso, and Nicolodi
\cite{JMN01}. For this we use the quaternionic model of M{\"o}bius
geometry \cite{BFLPP02,FLPP01, J03}.

\begin{figure}[h]
  \centering
  \scalebox{.4}{\includegraphics{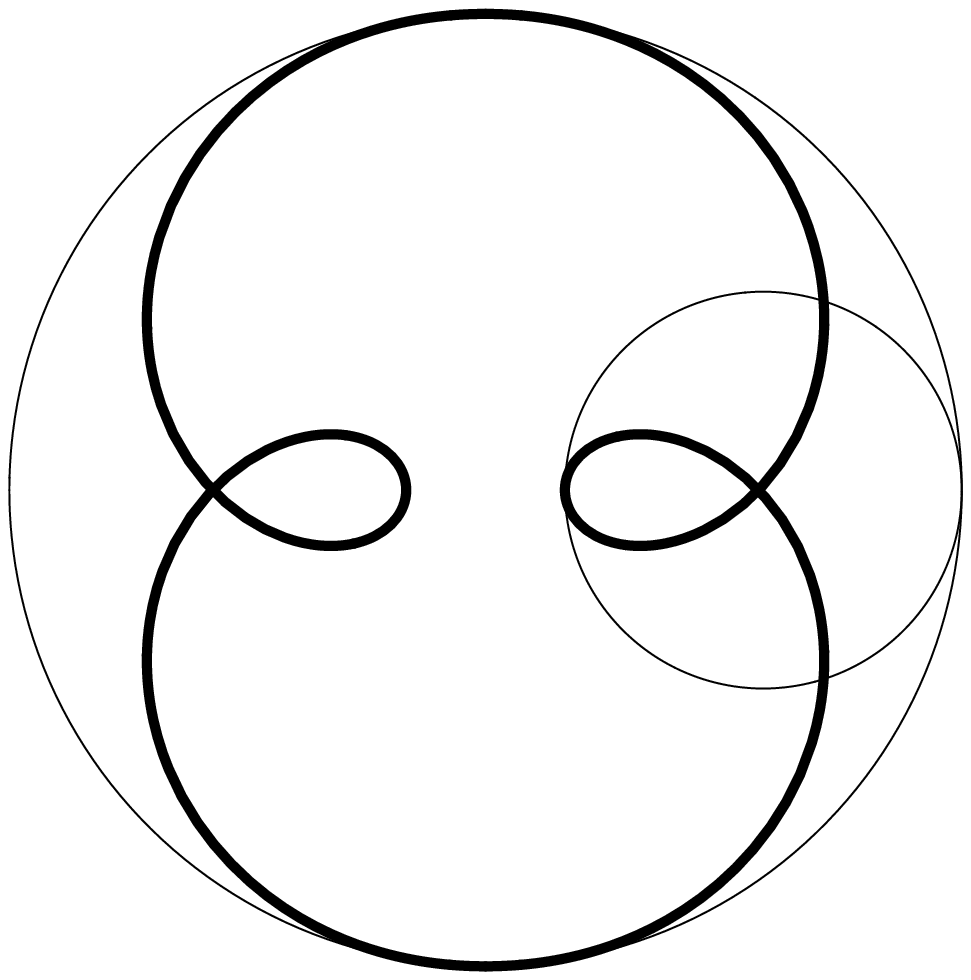}}\hfill
  \scalebox{.4}{\includegraphics{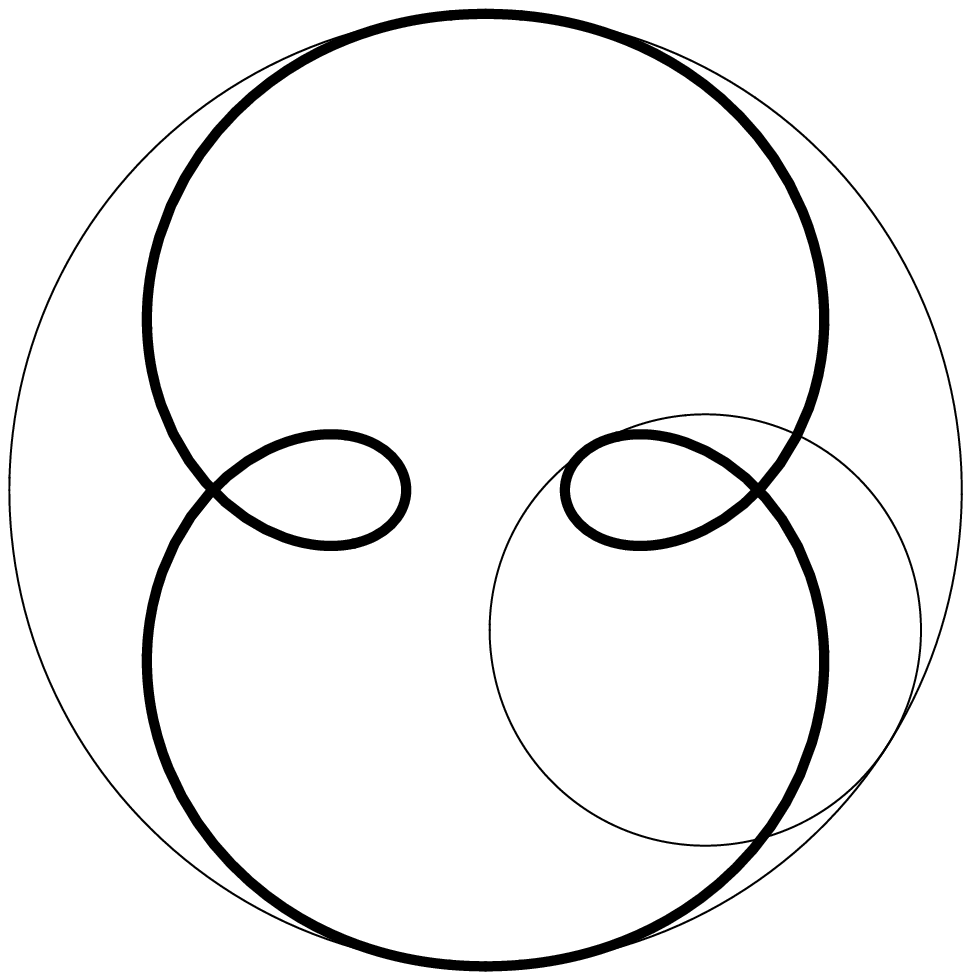}}\hfill
  \scalebox{.4}{\includegraphics{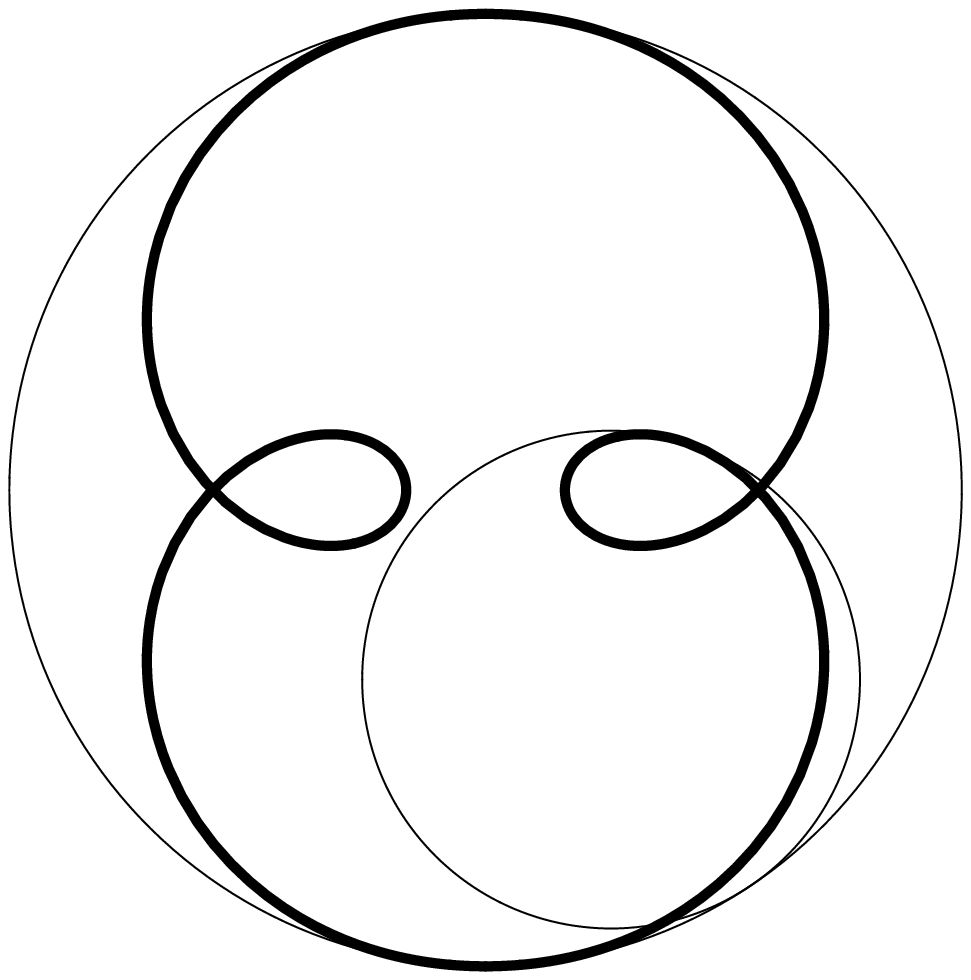}}
  \caption{Bryant surface with mean curvature spheres.}
  \label{fig:bryant_sphere_with_mcs}
\end{figure}

The main idea is a M{\"o}bius geometric characterization of
Bryant surfaces. A conformal immersion into $S^3$ is a constant mean
curvature $\pm 1$ surface in the hyperbolic space of curvature $-1$ with
asymptotic boundary $S^2\subset S^3$ if and only if all its mean curvature
spheres are tangent to $S^2$, i.e., all its mean curvature spheres are
horospheres, cf.\ Figure~\ref{fig:bryant_sphere_with_mcs}.  (This is
because horospheres are the spheres with mean curvature $\pm 1$ in
hyperbolic space.) The Bryant representation appears naturally in this
context when the hyperbolic Gauss map, which describes the intersection of
the mean curvature spheres with the asymptotic boundary, is interpreted
as a Darboux transform of the Bryant surface.

The quaternionic approach to 4--dimensional M{\"o}bius geometry, as introduced
in \cite{BFLPP02}, is based on the fact that the quaternionic projective line
$\HP^1$ with its standard conformal structure inherited from $\R^4=\H$ by the
decomposition
\begin{align*}
  \HP^1=\theset{\left[\lambda,1\right]}{\lambda\in\H}\cup\{[1,0]\}=\H\cup
  \{\infty\}
\end{align*}
is conformally equivalent to the standard 4--sphere and that, moreover, the
group of projective transformations corresponds to the group of orientation
preserving M{\"o}bius transformations. This is a 4--dimensional analogue to
the usual interpretation of $\CP^1=\C \cup \{\infty\}$ as the Riemann sphere.

In the following, we identify maps of a Riemann surface $M$ into $\HP^1$ with
line subbundles $L\subset\H^2$ of the trivial quaternionic rank 2 vector
bundle $\H^2$ over $M$. As explained in \cite[4.2]{BFLPP02}, such a map is a
conformal immersion if and only if $L$ is an \emph{immersed quaternionic
  holomorphic curve} in $\HP^1$, which means that the derivative $\delta=\pi
d_{|_L}\in\Omega^1\Hom(L,\H^2/L)$ of $L$ (with $\pi\colon \H^2\to\H^2/L$
the canonical projection) is nowhere vanishing and satisfies
\begin{equation}\label{eq:holomorphic_curve_condition}
  *\delta=\delta J
\end{equation}
for some $J\in\Gamma(\End(L))$ with $J^2=-1$ (where $*$ denotes precomposition
of the complex structure of $TM$). Slightly more general, a map $L$ from a
Riemann surface into $\HP^1$ is called a \emph{holomorphic curve} in $\HP^1$
if it admits $J$ with \eqref{eq:holomorphic_curve_condition}.

A fundamental object in M{\"o}bius geometric surface theory is the \emph{mean
  curvature sphere congruence} of a conformal immersion, which is the unique
congruence of oriented touching spheres with the conformally invariant
property that pointwise the immersion and the corresponding sphere have the
same mean curvature vector. In order to characterize the mean curvature sphere
congruence of a holomorphic curve in $\HP^1$ in the quaternionic language we
need the following description of 2--spheres \cite[3.4]{BFLPP02}: the oriented
2--spheres in $\HP^1$ are in one--to--one correspondence with the quaternionic
linear endomorphisms $S\in\End(\H^2)$ with $S^2=-\Id$. Such an endomorphism is
identified with the 2--sphere $\theset{[x]\in\HP^1}{[Sx]=[x]}$ which we denote
by $S$, too.  The endomorphism $S$ induces an orientation of the corresponding
2--sphere as it distinguishes a complex structure on the tangent bundle.
Hence, $S$ and $-S$ describe the same 2--sphere with different orientations.

As shown in \cite[5.2]{BFLPP02}, the mean curvature sphere congruence of $L$
is the unique section $S\in \Gamma(\End(\H^2))$ with $S^2=-\Id$ that satisfies
\begin{equation}\label{eq:mcs_condition}
  SL=L,\qquad *\delta=S\delta=\delta S,\qquad \textrm{and}\qquad Q_{|_L}=0
\end{equation}
where $Q=\tfrac14(SdS-*dS)$. The first two conditions describe oriented
touching of the immersion and the sphere at the corresponding point while the
third condition singles out the mean curvature sphere congruence among all
congruences of touching spheres. In the second equation, $S$ stands for the
induced complex structures on $\H^2/L$ and $L$, in particular, $S$ restricted
to $L$ equals the given complex structure $J$ of the holomorphic curve $L$
determined by \eqref{eq:holomorphic_curve_condition}.  Given the first two
conditions in \eqref{eq:mcs_condition}, the third one is equivalent to
$\im(A)\subset L$ where $A=\tfrac14(SdS+*dS)$.

Before we come to the M{\"o}bius geometric interpretation of Bryant's
representation we need to introduce the notion of Darboux transformation (in
the isothermic surface sense): two nowhere intersecting conformal immersions
of the same Riemann surface $M$ into the conformal 4--sphere are called a
\emph{Darboux pair} or \emph{Darboux transform} of each other if there is a
family of 2--spheres parametrized by $M$ that touches both immersions at
corresponding points with the right orientation. As we will see below
\eqref{eq:darboux_condition} this definition immediately generalizes to
holomorphic curves in $\HP^1$ that are not immersed.

Our treatment of Darboux transformations is in the spirit of \cite{BLPP,Boh03}
although we do not deal here with the more general notion of Darboux
transformation presented there. Cf.\ \cite{J03} for more information on the
Darboux transformation of isothermic surfaces.

Assume $L, L^\sharp\subset\H^2$ are two conformal immersions of a Riemann
surface $M$ into $\HP^1$ that do not intersect, i.e., $\H^2=L\oplus L^\sharp$.
Because $L^\sharp\cong\H^2/L$ and $L\cong\H^2/L^\sharp$ are canonically
isomorphic, it makes sense to interpret the derivatives $\delta$ of $L$ and
$\delta^\sharp$ of $L^\sharp$ as 1--forms with values in $\Hom(L,L^\sharp)$
respectively $\Hom(L^\sharp,L)$.  With respect to the splitting $\H^2=L\oplus
L^\sharp$, the trivial connection $d$ of $\H^2$ takes the form
\begin{align}\label{eq:splitt_connection}
  d=\dvector{\nabla^L&\delta^\sharp\\\delta&\nabla^\sharp}
\end{align}
where $\nabla^L$ and $\nabla^\sharp$ are connections induced on $L$ and
$L^\sharp$.  A sphere congruence $\tilde S$ that pointwise intersects both $L$
and $L^\sharp$ has to be of the form
\begin{equation}
  \label{eq:ribaucour}
  \tilde S=\dvector{J&0\\0&J^\sharp}
\end{equation}
with $J\in \Gamma(\End(L))$ and $J^\sharp\in \Gamma(\End(L^\sharp))$. The
condition that such a sphere congruence $\tilde S$ touches both $L$ and
$L^\sharp$ with the right orientation is (cf.\ second condition in
\eqref{eq:mcs_condition})
\begin{align}
  \label{eq:darboux_condition}
  {*}\delta &=J^\sharp \delta=\delta J& &\text{and}& 
  {*}\delta^\sharp &=J \delta^\sharp=\delta^\sharp J^\sharp. 
\end{align}
In particular, $J$ and $J^\sharp$ are the complex structures of the
holomorphic curves $L$ and $L^\sharp$. Obviously, \eqref{eq:darboux_condition}
makes sense for non--immersed holomorphic curves as well. In the following we
take \eqref{eq:darboux_condition} as the definition of Darboux transformations
in the context of (not necessarily immersed) holomorphic curves $L$,
$L^\sharp$ with $L\oplus L^\sharp =\H^2$.

An important characterization of the Darboux transformation is in terms of the
retraction form \cite{J03}: a \emph{retraction form} of a holomorphic curve
$L$ is a closed 1--form $\omega\in\Omega^1\End(\H^2)$ satisfying $\im
\omega\subset L\subset\ker\omega$.  For a retraction form $\omega$, the
connection $d-\omega$ is flat, because $d\omega=0$ and $\omega\land\omega=0$
imply the Maurer--Cartan equation $d\omega=\omega\land\omega$.  Moreover, at
the points where the mean curvature sphere $S$ of $L$ exists, e.g.\ where $L$
is immersed, $\omega$ satisfies
\begin{equation}
  \label{eq:retraction_complex}
  *\omega=S\omega=\omega S.
\end{equation}

Given two quaternionic holomorphic curves $L, L^\sharp\subset\H^2$ that form a
Darboux pair, a retraction form $\omega \in \Omega^1\End(\H^2)$ for $L$ can be,
with respect to the splitting $\H^2=L\oplus L^\sharp$, defined by
\begin{equation}
  \label{eq:def_of_tau}
  \omega =
  \begin{pmatrix}
    0 & \delta^\sharp \\ 0 & 0
  \end{pmatrix}.
\end{equation}
This retraction form $\omega$ has the property that $L^\sharp$ is a parallel
subbundle with respect to the flat connection $d-\omega$.  Conversely, given a
retraction form $\omega$ for $L$, every $(d-\omega)$--parallel line bundle in
$\H^2$ that does not intersect $L$ is a Darboux transform of $L$ (by
\eqref{eq:splitt_connection} and \eqref{eq:retraction_complex}).

We obtain that, for simply connected $M$, every Darboux transform $L^\sharp$
of a holomorphic curve $L$ is of the form.
\begin{equation}
  \label{eq:Darboux_with_F}
  L^\sharp =F\dvector{\k\\1}\H 
\end{equation}
where
$F\colon M\to\GL(2,\H)$ is a solution to 
\begin{equation}
  \label{eq:Fs_equation_tau}
  dF=\omega F
\end{equation}
for some retraction form $\omega \in \Omega^1\End(\H^2)$ of $L$.  The
Darboux transform $L^\sharp$ is immersed if and only if $\omega$ has
no zeros. The choice of the vector $\tvector{\k\\1}$ is arbitrary for
the present statement, but has the effect that
formula~\eqref{eq:quaternionic_halfspaceformula} coincides with
formula~\eqref{eq:Formula_in_halfspace}.  

For the rest of this section we will work with the 3--dimensional hyperbolic
space defined by fixing the following 2-- and 3--sphere in $\HP^1$: $\C \cup
\{\infty\} \subset \HP^1$ with corresponding endomorphism $S_0={\i}\Id$ and
$\Span_\R \{1,{\i},\k\} \cup \{\infty\} \subset \HP^1$, respectively.  The
hyperbolic space we are working with is thus the Poincar{\'e} half space
\begin{equation}
  \label{eq:half_space}
  \hyp=\theset{v_1+v_2{\i}+v_3\k}{v_1,v_2,v_3\in\R,\,v_3>0}\subset\H\cup
  \{\infty\}=\HP^1.  
\end{equation}
We call a surface in $\HP^1$ a \emph{Bryant surface} if it is, up to some
projective transformation, a Bryant surface in $\hyp$ in the usual sense. Note
that this includes surface with mean curvature minus one in $\hyp$.

The M{\"o}bius geometric description of the Bryant representation
in terms of Darboux transformations is given by the following theorem.

\begin{The}[\cite{JMN01}]\label{T:Bryant_representation}
  A non--totally umbilic surface immersed into $\HP^1$ is a Bryant surface if
  and only if it has a totally umbilic Darboux transform.  The totally umbilic
  Darboux transform is then the hyperbolic Gauss map of the Bryant surface.
\end{The}

Note that at a smooth Bryant end the immersion coincides with the hyperbolic
Gauss map which therefore ceases to be a Darboux transform, because the
splitting $L\oplus L^\sharp$ ``collapses''.

\begin{proof}
  Suppose that the conformal immersion $L^\sharp\subset\H^2$ has a totally
  umbilic Darboux transform $L\subset\H^2$. Because the notion of Darboux
  transformation is symmetric, $L^\sharp$ is also a Darboux transform of
  $L$. As explained above there is a retraction form $\omega$ of $L$ and a
  solution $F\colon \tilde M\to \GL(2,\H)$ to \eqref{eq:Fs_equation_tau}
  defined on the universal covering $\tilde M$ of $M$ such that
  $L^\sharp=F\tvector{\k\\1}\H$.  We may assume that $L$ is contained in the
  2--sphere $S_0$ and that for some point $p_0\in M$ we have $F_{p_0}=\Id$.
  The mean curvature sphere of $L$ is then $S_0$ and, by
  \eqref{eq:retraction_complex}, we have $*\omega=S_0\omega=\omega S_0$.  This
  shows that $\omega$ is a complex holomorphic 1--form with values in
  $\gl(2,\C)$.  Moreover, $\im(\omega)\subset L\subset\ker(\omega)$ implies
  that $\omega$ is traceless and has vanishing determinant, so $F$ is a
  holomorphic null immersion into $\SL(2,\C)$.
  
  The immersion $L^\sharp$ is either contained in $S_0$ or does not intersect
  $S_0$ at all, because for all $\lambda\in\H$ the section
  $\varphi:=S_0F\tvector{\k\\1}-F\tvector{\k\\1}\lambda$ is
  $(d-\omega)$--parallel and therefore vanishes identically or has no zeros at
  all. Since $L^\sharp$ is not totally umbilic, $L^\sharp$ and $S_0$ do not
  intersect.
  
  In order to prove that $L^\sharp$ is a Bryant surface, we have to show that
  the immersion $L^\sharp$ takes values in $\hyp$ and that the sphere
  congruence $\tilde S$ that touches $L$ and $L^\sharp$ (see
  \eqref{eq:ribaucour}) is the mean curvature sphere of $L^\sharp$. To do so
  we write the holomorphic null immersion $F$ into $\SL(2,\C)$ as
  \begin{equation}
    \label{eq:F_matrix}
    F=\dvector{a&b\\c&d}\colon \tilde M\to \SL(2,\C).
  \end{equation}
  The immersion  $L^\sharp=F\tvector{\k\\1}\H$ then becomes 
  \begin{equation}
    \label{eq:quaternionic_halfspaceformula}   
      L^\sharp=\frac{a\bar c+b\bar d+\k}{|c|^2+|d|^2}\colon M\to\hyp,
  \end{equation}
  which is exactly the formula~\eqref{eq:Formula_in_halfspace} for the Bryant
  representation in the half space model. Because $L^\sharp$ is defined on
  $M$, $F$ has monodromy in $\SU(2)$. This follows from the fact that a
  $\SL(2,\C)$--matrix is in $\SU(2)$ if and only if it leaves the quaternionic
  line $\tvector{\k\\1}\H$ invariant.

  To show that the 2--sphere congruence $\tilde S$, as given in
  \eqref{eq:ribaucour}, is the mean curvature sphere of $L^\sharp$ we have to
  check that
  \begin{equation}
    \label{eq:Q_of_ribaucour}
    Q=\frac14(\tilde Sd\tilde S-{*}d\tilde S)=\dvector{Q^L&0\\0&Q^\sharp}
  \end{equation}
  vanishes on $L^\sharp$, or, equivalently, that
  $Q^\sharp=\frac14(J^\sharp\nabla^\sharp J^\sharp-{*}\nabla^\sharp
  J^\sharp)=0$ where $\nabla^\sharp$ denotes the connection on
  $\End(L^\sharp)$ induced by the decomposition \eqref{eq:splitt_connection}
  of the trivial connection $d$. This fact is a consequence of
  $dS_0=0$: because $S_0$ is the mean curvature sphere of $L$, there is a
  section $H\in \Gamma(\Hom(L^\sharp,L))$ such that $S_0$, with respect to the
  splitting $\H^2=L\oplus L^\sharp$, takes the form
  \begin{equation*}
    S_0 =
    \begin{pmatrix}
      J & H \\ 0 & J^\sharp 
    \end{pmatrix}
    \qquad \textrm{ and therefore } 
    \qquad
    0=dS_0= 
    \begin{pmatrix}
      * & * \\ * & \nabla^\sharp J^\sharp + \delta H
    \end{pmatrix}.
  \end{equation*}
  Hence, $\nabla^\sharp J^\sharp =- \delta H$ and, by $*\delta = J^\sharp
  \delta$ (see \eqref{eq:darboux_condition}), we obtain $Q^\sharp=0$.
  
  To prove the converse suppose now that $L^\sharp$ is a Bryant
  surface in $\hyp$ and denote by $S$ the mean curvature sphere
  congruence of $L^\sharp$.  Then $S_p$ is a horosphere for all $p\in
  M$, i.e., it touches the ideal boundary $S_0$.  For every $p\in M$,
  define $L_p$ to be the point of intersection of $S_0$ and
  $S_p$. This defines a line bundle $L$ with $\H^2=L\oplus L^\sharp$
  and it suffices to check that $*\delta\psi=S\delta\psi=\delta S\psi$ for all
  $\psi\in\Gamma(L)$: then $L$ is a holomorphic curve with the complex
  structure induced by $S$, see
  \eqref{eq:holomorphic_curve_condition}, $L$ is totally umbilic,
  because its image is contained in $S_0$, and, by definition
  \eqref{eq:darboux_condition}, $L$ is a Darboux transform of
  $L^\sharp$ .
  
  Because both 2--spheres $S_p$ and $S_0$ have the same tangent space
  at $L_p$ we have ${S_p}_{|_{L_p}} ={S_0}_{|_{L_p}}$ and $S_p \equiv
  S_0 \mod L_p$ (we may assume that $S_0$ and $S_p$ touch with the
  same orientation).  Hence $R:=S-S_0$ satisfies $\ker R= L = \im R$,
  $L$ is smooth, and for $\psi \in\Gamma(L)$ we obtain
  \begin{equation}\label{eq:dS}
    dS(\psi)=d(S\psi)-Sd\psi=d(S_0\psi)-Sd\psi=(S_0-S)d\psi=R\delta\psi. 
  \end{equation}
  In particular $dS\psi\in\Omega^1(L)$. Because $dS=2({*}Q-{*}A)$ and $\im
  A\subset L^\sharp$ (with $Q$ and $A$ as in \eqref{eq:mcs_condition} and
  below) $dS\psi\in\Omega^1(L)$ implies $dS\psi =2{*}Q\psi$. In particular,
  $*dS\psi =-SdS\psi=dSS(\psi)$ and, using \eqref{eq:dS}, $*R\delta\psi
  =-SR\delta\psi=R\delta S(\psi)$.  Because $R$ is nowhere vanishing and
  anti--commutes with $S$, this proves $*\delta\psi=S\delta\psi=\delta S\psi$.
\end{proof}

\section{Soliton Spheres and \\ Quantization of the Willmore Energy}\label{sec:soliton_spheres}

In this section we prove the quantization of the Willmore energy for compact
Bryant surfaces with smooth ends and the fact that Bryant spheres with smooth
ends are soliton spheres.

As we have seen in Corollary~\ref{C:bryant_surface_with_smooth_ends}, a
compact Bryant surface with smooth ends is represented by a meromorphic
$\SL(2,\C)$--valued map $F$ with $\SU(2)$--monodromy on a compact Riemann
surface $M$. Although $F$ itself is not defined on $M$, but on its universal
covering $\tilde M$, the total pole order of $F$ on $M$ is well defined. The
following theorem shows that the Willmore energy of a compact Bryant surface
is quantized and directly related to this total pole order.

\begin{The}\label{T:total_curvature}
  The Willmore energy of a compact Bryant surface with smooth ends is
  $4\pi N$ with $N$ the total pole order of $F$ on a fundamental domain.
\end{The}

\begin{proof}
  We prove this theorem using the quaternionic Pl{\"u}cker
  formula~\cite{FLPP01} and thereby introduce the notation needed in the proof
  of Theorem~\ref{T:equality}. Let $L \subset \H^2$ be the holomorphic curve
  on a compact Riemann surface $M$ that, away from finitely many points,
  parametrizes the Bryant surface and denote by $L^\sharp\subset\H^2$ its
  hyperbolic Gauss map%
  \footnote{In contrast to the proof of Theorem~\ref{T:Bryant_representation}
    the roles of $L$ and $L^\sharp$ are interchanged here.}, %
  which, by Theorem~\ref{T:Bryant_representation}, is a Darboux transform of
  $L$.  By Corollary~\ref{C:bryant_surface_with_smooth_ends}, the Bryant
  representation $F$ of $L$ is a meromorphic $\SL(2,\C)$--valued map on
  $\tilde M$ with poles at the ends. In particular, the hyperbolic Gauss map
  $L^\sharp$ extends holomorphically through the ends, because it is the
  kernel of the meromorphic 1--form $\omega = dF F^{-1}$.  The fact that $F$
  has $\SU(2)$--monodromy implies that away from the ends
  \[
    \psi :=F\tvector{\k\\1}
  \]
  is a section of $L$ with quaternionic monodromy.

  We prove now that the section $\gamma$ of $L^{-1} =(\H^2)^*/L^\perp$ with
  monodromy defined by $\gamma(\psi)=1$ is holomorphic with respect to the
  unique quaternionic holomorphic structure $D$ on $L^{-1}$, see
  \cite{FLPP01}, for which the constant sections of $(\H^2)^*$ project to
  holomorphic sections of $L^{-1}$.  This holomorphic structure satisfies
  $D\pi =(\pi d)''$ where $\pi\colon (\H^2)^* \to L^{-1} =(\H^2)^*/L^\perp$
  denotes the canonical projection and $''$ denotes $\bar K$--part of the
  1--form. Let $\hat \gamma$ be the unique section of
  $(L^\sharp)^\perp\subset(\H^2)^*$ that away from the ends is defined by
  $\hat \gamma(\psi) =1$.  Then $\pi\hat\gamma=\gamma$ and
  $d\hat\gamma(\psi)=- \hat\gamma(d\psi)=-\hat\gamma(\omega\psi)=0$, since
  $\omega$, as a retraction form for $L^\sharp$, takes values in
  $L^\sharp$. Hence $d\hat\gamma$ takes values in $L^\perp$ and $D\gamma=D\pi
  \hat \gamma = (\pi d\hat \gamma)''=0$.  The section $\gamma$ extends
  smoothly through the ends and its vanishing order equals the pole order of
  $F$ (which can be seen using the normalization of Lemma~\ref{L:normal_form},
  because $\gamma=(\pi e_2^*) \bar z^n(\bar d - c\k )^{-1}$).
  
  Since $\gamma$ is a globally defined holomorphic section of $L^{-1}$ with
  monodromy, we can apply the quaternionic Pl{\"u}cker formula \cite{FLPP01}
  to the 1--dimensional linear system with monodromy spanned by $\gamma$: let
  $\nabla$ be the flat connection on $L$ defined away from the ends by
  $\nabla\psi=0$. Its $\bar K$--part $\nabla''$ is a quaternionic holomorphic
  structure on $L$ defined away from the ends, the $\bar K$--part of the dual
  connection on $L^{-1}$ coincides with the usual holomorphic structure $D$,
  because $\gamma$ is both $D$--holomorphic and $\nabla$--parallel.
  Therefore, the quaternionic Pl{\"u}cker formula implies that the difference
  of the Willmore energies of the quaternionic holomorphic line bundles
  $(L^{-1},D)$ and $(L,\nabla'')$ satisfies
  \begin{equation*}
    W(L^{-1},D)-W(L,\nabla'')=4\pi(-\deg(L^{-1}) + \ord (\gamma)).
  \end{equation*}
  Here, $\ord (\gamma)=N$ is the vanishing order of $\gamma$ on a fundamental
  domain.  Moreover, we have $\deg(L^{-1})=1-g$ with $g$ the genus of $M$,
  because $L$ is contained in some 3--sphere in $\HP^1$ and therefore $L$ is a
  quaternionic spin bundle \cite{PP98}.  Finally, $W(L,\nabla'')=0$, because
  the Hopf field $Q$ of $\nabla''$ on $L$ vanishes. (This fact has already
  been proven in the preceding section: $\nabla$ on $L$ as defined above
  coincides with the connection induced from $d$ by the splitting
  $\H^2=L\oplus L^\sharp$, see \eqref{eq:splitt_connection}, and because
  $L^\sharp$ is a totally umbilic Darboux transformation, the Hopf field of
  $\nabla''$ vanishes as proven in the paragraph following
  \eqref{eq:Q_of_ribaucour}.)  We obtain
  \begin{equation}\label{eq:willmore_energy_Linv}\begin{split}
      W(L^{-1})&=4\pi(g-1 + \ord (\gamma))=4\pi(g-1 + N).
    \end{split}\end{equation} 
  The Willmore energy of $L^{-1}$ satisfies $W(L^{-1})=\int_M(H^2-K)dA$ where
  $H$ is the mean curvature, K the Gaussian curvature, and $dA$ the area
  element with respect to the Euclidean geometry of
  $\R^3\cong\Span_\R\{1,{\i},\k\}\subset \H$, see~\cite{FLPP01}. Assuming
  that none of the ends lies at $\infty$ of $\R^3$,  the
  Gauss--Bonnet Theorem and \eqref{eq:willmore_energy_Linv} imply\\[-1ex]
  \begin{equation*}
    W=\int_MH^2dA =\int_M(H^2-K)dA+4\pi (1-g)=4\pi \ord(\gamma)=4\pi N. \\[-3ex]
  \end{equation*}
\end{proof}

For a Bryant sphere with smooth ends $F$ is rational, so
\eqref{eq:Formula_in_Ball} shows that Bryant spheres with smooth ends admit
conformal parametrizations in terms of rational functions. This is a
fundamental property of soliton spheres \cite{Pe04,BoPe}. A \emph{Soliton
  sphere} is a conformal immersion $L\subset\H^2$ of $\CP^1$ into $\HP^1$ such
that the linear system $(\H^2)^*\subset H^0(L^{-1})$ whose elements are the
homogeneous coordinates of $L$, the so called \emph{canonical linear system},
is contained in a linear system with equality in the quaternionic Pl{\"u}cker
estimate: for an $(n+1)$--dimensional linear system $H\subset H^0(L^{-1})$ of
a quaternionic holomorphic line bundle $L^{-1}$ over a compact Riemann surface
of genus $g$ the Pl{\"u}cker estimate \cite{FLPP01} states that the Willmore
energy of $L^{-1}$ satisfies
\begin{align*}
  \frac1{4\pi}W(L^{-1})\geq (n+1)(n(1-g)-\deg L^{-1}) +\ord(H).
\end{align*}
If $L$ lies in some 3--sphere in $\HP^1$, then $\deg L^{-1}=1-g$, cf.\
\cite{PP98}, and
\begin{align}\label{eq:pluecker_estimate_S3}
  \frac1{4\pi}W(L^{-1})\geq (n^2-1)(1-g)+\ord(H).
\end{align}
The above definition of soliton spheres is a generalization of the one
Iskander Taimanov, motivated by the soliton theory of the mKdV equation, gives
for spheres with special symmetry \cite{Ta99}.  In \cite{BoPe} we prove that
all Willmore spheres, i.e., all spheres obtained from complete minimal
surfaces of finite total curvature with planar ends, are soliton spheres. We
conclude the present article with a proof of the analogous result that all
Bryant spheres with smooth ends are soliton spheres. We actually prove a more
general statement for compact Bryant surfaces with smooth ends of arbitrary
genus.

\begin{The}\label{T:equality}
  The canonical linear system of a non--totally umbilic compact Bryant surface
  with smooth ends is contained in a 3--dimensional linear system with
  monodromy that has equality in the quaternionic Pl{\"u}cker estimate. In
  particular, Bryant spheres with smooth ends are soliton spheres.
\end{The}

\begin{proof}
  We proceed using the notation of the proof of
  Theorem~\ref{T:total_curvature}. Denote by $(\H^2)^*\subset H^0(L^{-1})$ the
  canonical linear system.  Define $\tilde H$ to be the 3--dimensional linear
  system with monodromy obtained by taking the span of $(\H^2)^*$ and the
  holomorphic section $\gamma$ defined in the proof of
  Theorem~\ref{T:total_curvature}.  The theorem is proven if we show that
  $\tilde H$ has equality in the quaternionic Pl{\"u}cker formula.
  
  Denote by $G\colon M\to S_0=\C\cup\{\infty\}\subset \H\cup\{\infty\}$, or,
  projectively $L^\sharp=\tvectork{G\\1}$, the hyperbolic Gauss map of $L$ and
  by $b(G)$ its total branching order.  The zeros of $\gamma$ are the ends of
  the Bryant surface and Lemmas~\ref{L:ord_at_regular_points} and
  \ref{L:ord_at_ends} below show that the Weierstrass order $\ord(\tilde H)$
  of $\tilde H$ is given by
  \begin{align}\label{eq:ord_tildeH}
    \ord(\tilde H)= \ord(\gamma)-2 \#\text{ends} + b(G).
  \end{align}
  
  The derivative $\delta^\sharp\in \Gamma(\Hom(L^\sharp,K\H^2/L^\sharp))$ is
  linear with respect to the complex structures induced by $S_0$ (see
  \eqref{eq:darboux_condition}) and it is complex holomorphic (see
  \cite{BFLPP02}).  Similarly, the 1--form $\omega=dF F^{-1}$ is a complex
  meromorphic section of $\Gamma(\Hom(\H^2/L^\sharp,KL^\sharp))$ (see the
  proof of Theorem~\ref{T:Bryant_representation}). Hence $\delta^\sharp
  \omega$ is a complex meromorphic quadratic differential on $M$, i.e., a
  meromorphic section of $K^2$. By definition, $\ord(\delta^\sharp)=b(G)$ and,
  because $\omega$ has second order poles at the ends (cf.\
  Lemma~\ref{L:poles_of_omega}) and no zeros, $\ord(\omega)=-2\#\text{ends}$.
  Using $\deg(K^2)=4(g-1)$, this implies
  \begin{align*}
    4(g-1)=\ord(\delta^\sharp\omega)=b(G)-2\#\text{ends}
  \end{align*}
  and, together with \eqref{eq:willmore_energy_Linv} and
  \eqref{eq:ord_tildeH}, we obtain 
  \begin{align*}
    \frac1{4\pi}W(L^{-1})=3(1-g)+\ord(\tilde H)
    =(n^2-1)(1-g)+\ord(\tilde H)
  \end{align*}
  where $n+1=\dim\tilde H=3$. Thus for $\tilde H$ equality holds in
  \eqref{eq:pluecker_estimate_S3}.
\end{proof}

The following two lemmas are needed to compute the total order
\eqref{eq:ord_tildeH} of the linear system $\tilde H$.

\begin{Lem}\label{L:ord_at_regular_points}
  If $p\in M$ is not an end, then $\ord_p(\tilde H)= b_p(G)$.
\end{Lem}

\begin{proof}
  In the proof of Theorem~\ref{T:total_curvature} we have, away from the ends,
  defined the section $\hat\gamma$ of $(L^\sharp)^\perp$. It satisfies $\pi
  \hat \gamma= \gamma$ and $\pi d \hat \gamma=0$ for the canonical projection
  $\pi\colon(\H^2)^*\to(\H^2)^*/L^\perp$. This means that $\hat\gamma$ can be
  interpreted as the unique prolongation of the holomorphic section $\gamma$
  to the 1--jet bundle of $L^{-1}$, cf.\ \cite{FLPP01}.
  
  We chose $\hat \alpha\in L^\perp_p$ and $\hat \beta=\hat \gamma_p\in
  (L^\sharp)^\perp_p$. Then $\alpha=\pi \hat \alpha$ and $\beta=\pi \hat
  \beta$ form a basis of the canonical linear system $(\H^2)^*$.  Because
  $\beta_p=\gamma_p$, the section $\beta$ does not vanish at $p$ and there are
  quaternion valued functions $f$ and $h$ defined in a neighborhood of $p$
  such that $\alpha=\beta f$ and $\gamma=\beta h$. The sections $\alpha$ and
  $\tilde \gamma = \gamma-\beta= \beta(h-1)$ both vanish at $p$. The section
  $\alpha$ vanishes to first order at $p$, since $\alpha$ and $\beta$ form a
  basis of the linear system $(\H^2)^*$, which has no Weierstrass points,
  because $L$ is immersed. In particular, $d_pf \neq 0$.
  
  In order to relate the vanishing order of $\tilde \gamma$ at $p$ to the
  branching order of the hyperbolic Gauss map $G$ at $p$, we use that
  $\hat\gamma = \hat \beta h + (\hat \alpha-\hat\beta f) \lambda$ for a
  quaternion valued function $\lambda$, since $\hat \alpha-\hat\beta
  f\in\Gamma(L^\perp)$. By definition of $\hat \beta$, we have $\lambda_p=0$.
  Since $d \hat \gamma$ takes values in $L^\perp$, the equation
  \begin{equation*}
    d\hat \gamma= \hat \beta dh -\hat \beta df \lambda+
    (\hat \alpha-\hat\beta f) d\lambda
  \end{equation*}
  implies $dh=df\lambda$ and $d\hat \gamma=(\hat \alpha-\hat\beta f)
  d\lambda$. For the vanishing order of $\tilde \gamma = \beta(h-1)$ at $p$ we
  therefore obtain (using $dh=df\lambda$, $d_pf\neq 0$ and $\lambda_p=0$)
  \begin{equation*}
    \ord_p(\tilde \gamma) = \ord_p(dh)+1 = \ord_p(\lambda) + 1 =
    \ord_p(d\lambda) + 2.
  \end{equation*}
  For the branching order of $L^\sharp$ we obtain, using $d\hat \gamma=(\hat
  \alpha-\hat\beta f) d\lambda$ and $\hat
  \gamma_p\neq0$, 
  \begin{equation*}
    b_p(G) = \ord_p(d\hat\gamma)= \ord_p(d\lambda).
  \end{equation*}
  The Weierstrass gap sequence of the linear system $\tilde H$ for
  the point $p$ is $0$, $1$, $b_p(G)+2$ (realized by the sections $\beta$,
  $\alpha$, $\tilde \gamma$). 
\end{proof}

\begin{Lem}\label{L:ord_at_ends}
  If $p\in M$ is an end, then $\ord_p(\tilde H)= \ord_p(\gamma)-2+b_p(G)$.
\end{Lem}

\begin{proof}
  Let $z\colon U\to\Delta$ be a holomorphic coordinate on an open neighborhood
  $U\subset M$ of $p$ such that $z(p)=0$.  Let $F=z^{-n}\tvector{a&b\\c&d}$ be
  the holomorphic null immersion representing the Bryant surface
  $L_{|_{\Delta^*}}$ and suppose that $F$ has the normal form of
  Lemma~\ref{L:normal_form}. Then $\ord_0(\gamma)=n\neq0$.  The hyperbolic
  Gauss map $L^\sharp\subset\H^2$ is given by the image of
  $\omega=F'F^{-1}dz$, see Section~\ref{sec:darboux}. Thus
  \begin{equation*}
    L^\sharp=\tvector{z^{-n-1}(zb'-nb)\\ z^{-n-1}(zd'-nd)}\H
    =\tvector{zb'-nb\\ zd'-nd}\H.
  \end{equation*}
  Since $\ord_0( b)=1$ (see the proof of Lemma~\ref{L:poles_of_omega}) and
  $d(0)\neq0$, it follows that either $n\geq 2$ and $b_0(G)=0$ or $n=1$ and
  $b_0(G)=\ord_0(b-b'(0)z)-1\geq 1$. In the first case the claimed formula
  follows, since $L$ is immersed and therefore $(\H^2)^*$ is Weierstrass point
  free.
  
  To see the formula in the case $n=1$, let $\alpha,\beta\in (\H^2)^*\subset
  H^0(L^{-1})$ be the projections to the first and second coordinate of
  $\H^2$. The section $\beta$ does not vanish at $z=0$, because $F$ is in
  normal form and therefore
  $L=[\psi]=\tvectork{z^{-1}a\k+z^{-1}b\\z^{-1}c\k+z^{-1}d}$ is
  $\tvectork{0\\1}$ for $z=0$.  The section $\gamma$ vanishes to order $1$.
  We now compute the vanishing order of $\alpha-\gamma \overline{b'(0)}$. We
  have
  \begin{equation*}    
    \alpha = \beta \bar z(\overline{c\k +d})^{-1} (\overline{a\k +b})\bar
    z^{-1} \qquad \textrm{and} \qquad 
    \gamma = \beta \bar z(\overline{c\k +d})^{-1}
  \end{equation*}
  where the last formula follows from
  \[\gamma\tvector{z^{-1}a\k+z^{-1}b\\z^{-1}c\k+z^{-1}d}=\gamma(\psi)=1.\] Thus
  $\alpha-\gamma \overline{b'(0)}=\beta \bar z(\overline{c\k
    +d})^{-1}(\overline{a\k +b-b'(0)z})\bar z^{-1}$ and, by $d(0)\neq0$, we
  have $\ord_0(\alpha-\gamma \overline{b'(0)})=\min(\ord_0 (a),\ord_0
  (b-b'(0)z))$. By Lemma~\ref{L:normal_form} (and its proof) $a'(0)=0$ and
  $a'd'-b'c'=1$, which shows $b'(0)c'(0)=-1$.  Using $b'(0)c'(0)=-1$, the
  identities $a'd'-b'c'=1$ and $ad-bc=z^2$ imply $\ord_0 (a)\geq \ord_0
  (b-b'(0)z)$.  This proves the claim, because then $\ord_0(\alpha-\gamma
  \overline{b'(0)})=\ord_0 (b-b'(0)z)=b_0(G)+1$. The Weierstrass order at $p$
  is therefore $\ord(\tilde H )_p = b_0(G)-1$, since the Weierstrass gap
  sequence at zero is $0,1,b_0(G)+1$ (realized by $\beta,\gamma,\alpha-\gamma
  \overline{b'(0)}$).
\end{proof}

\begin{Rem}
  At the beginning of the proof of Lemma~\ref{L:ord_at_ends} we have seen that
  the hyperbolic Gauss map is immersed at a smooth end if and only if the pole
  order of $F$ at the end is greater than $1$, i.e., if the end is catenoidal
  and not horospherical in the sense of \cite{ET01}.
\end{Rem}

\section{Bryant's quartic differential}\label{sec:quartic_differential}

We show that the vanishing of Bryant's quartic differential~$\mathcal Q$
provides a uniform Möbius geometric characterization of the compact immersed
surfaces $f\colon M\to S^3$ in the conformal 3--sphere~$S^3$ that can be
obtained as compactification of either a Euclidean minimal surface with planar
ends or a Bryant surfaces with smooth ends.

The quartic differential $\mathcal Q$ was introduced by Bryant \cite{Br84} for
Willmore surfaces in the conformal 3--sphere $S^3$.  It was pointed out by
Konrad Voss in a talk given at Oberwolfach that $\mathcal Q$ may be defined
for arbitrary conformal immersions of a Riemann surface $M$ into $S^3$ and
that $\mathcal Q$ is holomorphic if and only if, locally and away from
umbilics and isolated points, the immersion is Willmore or has constant mean
curvature with respect to some space form geometry, see \cite{Boh08} for a
proof. We prove the following global characterization of compact surfaces with
$\mathcal Q\equiv0$.

\begin{The}\label{T:vanishing_of_Q} Let $f\colon M\to S^3$ be a conformal
  immersion of a compact Riemann surface $M$ into the conformal
  3--sphere~$S^3$. The quartic differential $\mathcal Q$ of $f$ vanishes
  identically if and only if $f$ is the compactification of either a Euclidean
  minimal surface with planar ends or a Bryant surface with smooth ends.
\end{The}

Here an immersion $f\colon M\to S^3$ into the conformal 3--sphere is called a
Euclidean minimal surface if and only if it admits a point $\infty\in S^3$ not
on $f$ in which all its mean curvature spheres intersect, i.e., if the
resulting immersion into $\R^3=S^3\backslash \{\infty\}$ is Euclidean
minimal. An immersion is the compactification of a Euclidean minimal surface
with planar ends if and only if the resulting immersion into
$\R^3=S^3\backslash \{\infty\}$ is a complete Euclidean minimal surface with
finite total curvature and planar ends, see \cite{Br84}.  Similarly, an
immersion $f\colon M\to S^3$ into the conformal 3--sphere is a Bryant
surface if and only if it has a totally umbilic darboux transform, 
cf.~Theorem~\ref{T:Bryant_representation}.

\begin{proof} 
  Lemma~\ref{L:umbilics} below implies that $f$ has only finitely many
  umbilics. The theorem then follows from
  Proposition~\ref{P:local_vanishing_of_Q}.
\end{proof}

Using that a holomorphic quartic differential on the sphere vanishes
identically and that, by \cite{Br84}, every Willmore spheres is the
compactification of a Euclidean minimal surface with planar ends, we obtain:

\begin{Cor}Let $f\colon \CP^1\to S^3$ be a conformal immersion of the
  2--sphere into the 3--sphere. The quartic differential $\mathcal Q$ of $f$
  is holomorphic if and only if $f$ is a Willmore sphere or a Bryant sphere
  with smooth ends.
\end{Cor}

For the proof of Lemma~\ref{L:umbilics} and
Proposition~\ref{P:local_vanishing_of_Q} we use the light cone model of the
conformal 3--sphere and the invariants introduced in \cite{BPP02}.  The light
cone model of the conformal 3--sphere $S^3$ is based on the identification
\begin{displaymath}
	S^3\cong \mathbb P(\mathcal L),\quad x\leftrightarrow [1:x] 
\end{displaymath}
of $S^3\subset \R^4$ with the projectivized light cone $\mathcal
L\subset\R^{4,1}$ in 5--dimensional Minkowski space $\R^{4,1}$ with the metric
$\langle x, x\rangle=-x_0^2+\sum_{i=1}^4 x_i^2$. Under this identification,
the group of Möbius transformations of $S^3$ is identified with the identity
component of the group of linear transformation that preserve the Minkowski
product $\langle\,,\rangle$.

Let $f\colon M\to S^3$ be a conformal immersion of a Riemann surface $M$ into
the conformal 3--sphere and $z\colon M\supset U\to \C$ a holomorphic chart
on~$M$. Then there is a unique holomorphic lift $\psi\colon U\to \mathcal L$
of $f$ such that
\begin{displaymath}
	\langle d\psi,d\psi\rangle =|dz|^2
\end{displaymath}
and $\psi_0>0$. Denote $\gamma\colon M\to \R^{4,1}$ the unique (up to sign)
smooth map with
\begin{displaymath}
	\gamma \perp \{\psi,\psi_z,\psi_{\bar z},\psi_{z\bar z}\}
	\text{\qquad and \qquad}
	\langle \gamma,\gamma\rangle = 1. 
\end{displaymath}
Then $\mathbb P(\gamma^\perp\cap \mathcal L)$ is the mean curvature sphere
congruence or conformal Gauss map of $f$.  The \emph{conformal Hopf
  differential} and \emph{Schwarzian derivative} of $f$ with respect to $z$
are the smooth functions $\kappa,c\colon U\to \C$ satisfying
\begin{equation}
	\psi_{z z}+\tfrac c2\psi =\kappa \gamma.
\end{equation}
The corresponding functions $\tilde \kappa,\tilde c$ with respect to another
holomorphic chart $\tilde z\colon U\to\C$ are given by
\begin{equation}\label{eq:transformation_rules}
	\tilde \kappa \frac{d\tilde z^2}{|d\tilde z|}=\kappa \frac{dz^2}{|dz|}
	\text{\qquad and\qquad}
	\tilde c d\tilde z^2= (c-S_z(\tilde z))dz^2, 
\end{equation}
where $S_z(\tilde z)=\left(\frac{\tilde z_{z z}}{\tilde
    z_z}\right)_z-\frac12\left(\frac{\tilde z_{zz}}{\tilde z_z}\right)^2$ is
the usual Schwarzian derivative.

Denote $\hat\psi\colon U\to\mathcal L$ the unique map such that the Minkowski
product takes the form
\begin{displaymath}
	\left(\begin{smallmatrix}
    	0& 0& 0& -1& 0\\
    	0& 0& \frac12& 0& 0\\
    	0& \frac12& 0& 0& 0\\
    	-1& 0& 0& 0& 0\\
    	0& 0& 0& 0& 1
    \end{smallmatrix}\right)
\end{displaymath}
with respect to the frame $(\psi,\psi_z,\psi_{\bar z},\hat\psi,\gamma)$.  The
frame equations are then
\begin{equation}\label{eq:frame_equations}\begin{aligned}
    \psi_{zz} &= -\tfrac c2\psi +\kappa\gamma,&
    \psi_{z\bar z} &= - |\kappa|^2\psi +\tfrac12\hat\psi,\\
    \hat\psi_z &= -2|\kappa|^2\psi_z-c\psi_{\bar z} + 2\kappa_{\bar z}\gamma,&
    \gamma_{z} &= 2\kappa_{\bar z}\psi- 2 \kappa \psi_{\bar z}
\end{aligned}\end{equation}
with compatibility conditions
\begin{align}
  &\qquad \tfrac12 c_{\bar z} = 2(|\kappa|^2)_z+2\bar\kappa_z\kappa
  &&\text{(Gauss equation),}\label{eq:gauss} \\
  &\qquad \Im(\kappa_{\bar z\bar z}+\tfrac{\bar c}2\kappa)=0&&\text{(Codazzi
    equation).}
\end{align}
In particular, the second $z$--derivative of $\gamma$ is
\begin{equation}
  \gamma_{z z}=(2\kappa_{z\bar z} +2\kappa|\kappa|^2)\psi + 2\kappa_{\bar
    z}\psi_z-2\kappa_z\psi_{\bar z} -\kappa\hat\psi. 
\end{equation}
\emph{Bryant's quartic differential} is defined as
\begin{equation}\label{eq:q}
  \mathcal Q=\tfrac14\langle \gamma_{z z},\gamma_{z z}\rangle dz^4
  = (\kappa\kappa_{\bar z z}+\kappa^2|\kappa|^2-\kappa_{\bar z}\kappa_z)dz^4
\end{equation}
(since $(\psi,\psi_x,\psi_y,\gamma,\hat\psi)$ with $z=x+iy$ is a local section
of $\mathcal F_f^{(\gamma)}$ in \cite{Br84}, Theorem~B).

The umbilic points of $f$ are those points at which the mean curvature sphere
touches $f$ to second order, i.e., umbilic points are the zeros of $\kappa$.
On the complement $M_0$ of the set of umbilic points of $f$, its mean
curvature sphere congruence has a unique \emph{second envelope}
$f^\sharp\colon M_0\to S^3$ which is defined by the property that any lift
$\psi^\sharp$ of $f^\sharp$ satisfies
\begin{equation}\label{eq:second_envelope_defining_equation}
  \langle\psi^\sharp,\psi^\sharp\rangle=\langle \psi^\sharp,\gamma\rangle
  = \langle \psi^\sharp,d\gamma \rangle = 0.
\end{equation}  
With respect to a chart $z$, a lift of $f^\sharp$ can be obtained by
\begin{equation}\label{eq:second_envelope}
  \psi^\sharp = 2|\kappa_{\bar z}|^2\psi - 2 \kappa_{\bar z}\bar\kappa\psi_z 
  -2\bar\kappa_z\kappa\psi_{\bar z}+ |\kappa|^2\hat\psi.
\end{equation}

\begin{Lem}\label{L:umbilics} 
  A conformal immersion $f\colon M\to S^3$ of a compact Riemann surface $M$
  into the 3--sphere whose quartic differential $\mathcal Q$ vanishes
  identically is either totally umbilic or has only finitely many umbilic
  points.
\end{Lem}

\begin{proof}
  The idea of the proof is to define a complete metric $g$ on the complement
  $M_0$ of the set of umbilics that has non--positive curvature and finite
  total curvature.  From~\cite[Chapter~III, Proposition~16]{La80} we then
  obtain that $M\setminus M_0$ consists of a finite number of points. In order
  to define the metric $g$ we use the Uniformization Theorem and distinguish
  the cases that $M$ has genus 0, 1, or greater than 1.
  \begin{enumerate}[1.]
  \item   Case $M=\CP^1$: without loss of generality we can assume that $\infty
  \in\CP^1$ is an umbilic point, i.e., $\infty \not\in M_0$. Denote $\kappa$
  the conformal Hopf differential with respect to a holomorphic chart $z\colon
  \CP^1\setminus\{\infty\} \to \C$. Then
  \begin{equation*}
  	g=e^{2u}|dz|^2,\qquad u=-\log|\kappa|
  \end{equation*}
  defines a (coordinate dependent) metric on $M_0$.  By
  \eqref{eq:transformation_rules} we have $g=\frac{|d\tilde
    z|^2}{|\tilde\kappa|^2|\tilde z|^8}$ for $\tilde z=\frac1z$. Thus $g$ is
  complete, because $\tilde \kappa$ is smooth near $\infty$, $\kappa$ is
  smooth on $\CP^1\setminus\{\infty\}$, and zero on $\CP^1\setminus M_0$. The
  curvature of $g$ is
  \begin{align*}
    K = - 4 e^{-2u} u_{z\bar z} = 4 |\kappa|^{2} \Re\left(\frac{\kappa_{\bar z
        z}\kappa -\kappa_z\kappa_{\bar z}}{\kappa^2}\right) =
    -4|\kappa|^4,
  \end{align*}
  where the last equality holds because $\mathcal Q\equiv 0$ in
  \eqref{eq:q}. The total curvature of $g$ is finite, because
  \begin{equation*}
  	\int_{M_0} KdA=-4\int_\C|\kappa|^2dx\wedge dy = -4 \tilde W,
  \end{equation*}
  where $W=\tilde W+2\pi\chi(M)$ is the Willmore energy of $f$, 
  cf. \cite[Section~3.2]{BPP02} and Appendix~\ref{sec:appendix}. 
  
  \item Case $M=\C/\Gamma$: the conformal Hopf differential $\kappa$ with
  respect to the chart $z$ on the universal cover $\C$ is invariant under the
  translations by elements in the lattice $\Gamma$ and defines a function
  $\kappa\colon M\to\C$. Hence as in the genus 0 case
  \begin{equation*}
  	g=e^{2u}|dz|^2,\qquad u=-\log|\kappa|
  \end{equation*}
  is a metric on $M_0$ with the desired properties.

  \item Case $M=\mathbf B^2/\Gamma$ for $\mathbf B^2\subset\C$ the unit disc
  with the Poincar\'e metric $\hat g=\frac{4|dz|^2}{(1-|z|^2)^2}$ and $\Gamma$
  a discrete group of isometries of $\mathbf B^2$: now $\kappa$ with respect
  to the chart $z$ on the universal cover $\mathbf B^2$ is not invariant under
  $\Gamma$, but
  \begin{equation*}
  	|\kappa|(1-|z|^2)\colon M\to\R
  \end{equation*}
  is well defined by \eqref{eq:transformation_rules}, because $\Gamma$ is a
  group of isometries of $\hat g$. Thus
  \begin{equation*}
  	g=e^{2u}|dz|^2=\frac{\hat g}{|\kappa|^2(1-|z|^2)^2},\qquad
  	u=-\log\left|\tfrac12\kappa(1-|z|^2)^2\right|
  \end{equation*}
  defines a complete metric on $M_0$ with non--positive curvature
  \begin{align*}
	K &= - 4 e^{-2u} u_{z\bar
          z}=-|\kappa|^4(1-|z|^2)^4-2|\kappa|^2(1-|z|^2)^2. 
  \end{align*}
  The total curvature of $g$ is finite, because
  \begin{align*}
    \int_{M_0} KdA &=-4\int_\Delta|\kappa|^2dx\wedge dy -2\int_\Delta\frac
    4{(1-|z|^2)^{2}}dx\wedge dy =-4\tilde W-2 \hat A,
  \end{align*}
  where $\Delta\subset \mathbf B^2$ is a fundamental domain of $\Gamma$ and
  $\hat A$ its hyperbolic area.
  \end{enumerate}
\end{proof}

The following proposition is a local version of Theorem~\ref{T:vanishing_of_Q}
which holds away from umbilic points.

\begin{Pro}\label{P:local_vanishing_of_Q} 
  Let $f\colon M\to S^3$ be an umbilic free conformal immersion of a Riemann
  surface $M$ into the conformal 3--sphere $S^3$. The quartic differential
  $\mathcal Q$ of $f$ vanishes identically if and only if either $f$ is
  Euclidean minimal in $\R^3=S^3\backslash \{\infty\}$ for some $\infty \in
  S^3$ or $f$ is a Bryant surface.
\end{Pro}

\begin{proof}
  By \eqref{eq:frame_equations} and \eqref{eq:q}, the derivative of
  $\psi^\sharp$ in \eqref{eq:second_envelope} is
  \begin{gather}\label{eq:psi_sharp_z}
    \psi^\sharp_z=\kappa^{-1}\kappa_z\psi^\sharp+\lambda\bar\kappa\gamma_z +
    \kappa^{-1}q\gamma_{\bar z},\qquad \intertext{ where }
    \lambda=\bar\kappa^{-1}\bar\kappa_{zz}+\tfrac c2\qquad \text{and}\qquad
    \mathcal Q= q \, dz^{4}.\notag
  \end{gather}
  If $\mathcal Q\equiv0$ then $\lambda$ vanishes at a point if and only if the
  second envelope $f^\sharp$ of the mean curvature sphere congruence of $f$ is
  not immersed at this point. The latter property is clearly independent of
  the chart. Because
  \begin{equation*}
  	\lambda_{\bar z}
  	=\bar\kappa^{-2}(\bar\kappa\bar\kappa_{\bar
  	zz}+\bar\kappa^2|\kappa|^2-\bar\kappa_{\bar z}\bar\kappa_z)_z
  	=\bar\kappa^{-2}(\bar q)_z=0
  \end{equation*}
  the function $\lambda$ is holomorphic and hence either vanishes identically
  or has isolated zeros. So  $f^\sharp$ is either constant or, away from
  isolated points, an immersion.  In the first case $f$ is Euclidean minimal
  in $\R^3=S^3\backslash\{\infty\}$ for $\infty = f^\sharp$ (note that $f$
  does not go through $\infty$, see
  \eqref{eq:second_envelope}).
  
  In case $f^\sharp$ is non--constant we assume for a moment that it is
  globally immersed. Then $f^\sharp$ is conformal, because $\langle
  \psi^\sharp_z,\psi^\sharp_z\rangle=0$ by
  \eqref{eq:second_envelope_defining_equation} and \eqref{eq:psi_sharp_z}, and
  $\langle \gamma_z,\gamma_z\rangle=0$ which follows from
  \eqref{eq:frame_equations}.  Taking the derivative in \eqref{eq:psi_sharp_z}
  shows that $\psi^\sharp_{zz}$ lies in the span of $\psi^\sharp$, $\gamma_z$,
  and $\gamma_{zz}$. Since $\mathcal Q\equiv0$, equations \eqref{eq:q} and
  \eqref{eq:second_envelope_defining_equation} imply that the complex subspace
  spanned by $\psi^\sharp$, $\gamma_z$, and $\gamma_{zz}$ is a null space for
  the non--degenerate symmetric product induced by the Minkowski metric.  But
  a null space for a non--degenerate symmetric product on a 5--dimensional
  space is at most 2--dimensional such that $\psi^\sharp$, $\psi^\sharp_z$,
  and $\psi^\sharp_{zz}$ are linearly dependent and $f^\sharp$ is totally
  umbilic. In other words, the second envelope $f^\sharp$ is a holomorphic map
  into a round 2--sphere in $S^3$. By continuity, this more generally holds in
  case $f^\sharp$ is globally immersed except for isolated points.

  Using \eqref{eq:psi_sharp_z} one can check that away from the branch points
  of $f^\sharp$ the frame $(\psi,\psi_x, \psi_y,\gamma,\hat\psi)$ induces the
  same orientation as $(\psi^\sharp,\psi_x^\sharp,
  \psi_y^\sharp,\gamma,\hat\psi)$, i.e., the envelopes $f^\sharp$ and $f$ of
  the mean curvature sphere congruence of $f$ touch with the same
  orientation. Hence $f^\sharp$ is a totally umbilic Darboux transform of $f$
  and $f$ is a Bryant surface by Theorem~\ref{T:Bryant_representation}.
  
  Conversely, if $f$ is Euclidean minimal then $f^\sharp$ is constant and
  \eqref{eq:psi_sharp_z} implies $\mathcal Q\equiv0$. If $f$ is a Bryant
  surface then Theorem~\ref{T:Bryant_representation} implies that $f^\sharp$
  is a totally umbilic branched conformal immersion, such that
  \begin{equation*}
  	0=\langle \psi^\sharp_z,\psi^\sharp_z\rangle=4\bar\kappa^2\lambda q,
  \end{equation*}
  by \eqref{eq:psi_sharp_z} and \eqref{eq:frame_equations}.  Hence $q$
  vanishes identically: otherwise there would exist an open set on which
  $\lambda$ vanished identically and $q$ had no zeros. But this and
  \eqref{eq:psi_sharp_z} would imply that $\psi^\sharp_{zz}$ is linearly
  independent of $\psi^\sharp$ and $\psi^\sharp_z$, which contradicts the
  assumption that $f^\sharp$ is totally umbilic.
\end{proof}

\begin{Rem}
  Let $f\colon M \rightarrow S^3$ be an immersion obtained from a Euclidean
  minimal surface with planar ends or a Bryant surface with smooth ends by
  filling in points at the ends. Then all added points are umbilic points,
  because away from umbilic points the second envelope $f^\sharp$ of the mean
  curvature sphere congruence of $f$ does never intersect $f$, see
  \eqref{eq:second_envelope}: since $\mathcal{Q}\equiv 0$, either $f^\sharp$
  is a constant point $\infty$ (in the Euclidean minimal case) or $f^\sharp$
  is the restriction of the hyperbolic Gauss map to the complement of the set
  of umbilic points (in the Bryant case). But in the Euclidean minimal case
  the planar ends of $f$ are the points where $f$ goes through $\infty$ and in
  the Bryant case the smooth Bryant ends of $f$ are the points at which $f$
  coincides with the hyperbolic Gauss map (see the remark following
  Theorem~\ref{T:Bryant_representation}).
\end{Rem}

The Willmore energy of an immersion $f\colon M \rightarrow S^3$ with
$\mathcal{Q}\equiv 0$ of a compact surface into the conformal 3--sphere is
always an integer multiple of $4 \pi$ (it is essentially $4\pi$ times the
number of ends, except that smooth catenoid cousin ends have to be counted
with the pole order of $F$, see Theorem~\ref{T:total_curvature}).  If
$M=\CP^1$, the possible Willmore energies are $W=4\pi d$, $d\in \N \backslash
\{0,2,3,5,7\}$, see \cite{Br84,Br87} for the Euclidean minimal case and
Theorem~\ref{T:gaps_for_spheres} for the Bryant case.

\appendix
\section{}\label{sec:appendix}

The total absolute curvature plays an important role in the theory of
Euclidean  minimal and Bryant surfaces. In this appendix we recall its relation
to the Willmore energy and show that if a Bryant
surface in $\ball$ is the intersection of $\ball$ with a compact surface in
$\R^3$ then it is a compact Bryant surface with smooth ends. 

Consider an immersion $f\colon M\to\bar M$ of an oriented
2--dimensional manifold $M$ into a 3--dimensional Riemannian manifold
$\bar M$. Let $dA$ be the area element, $H$ the mean curvature, $G$
the Gaussian curvature (the determinant of the Weingarten operator), $K$
the curvature of the induced metric on $M$, and $\bar K$ the sectional
curvature of $\bar M$ on the tangent spaces of $f$. The latter three
quantities are related by the Gauss equation $K=G+\bar K$. The 1--form
\begin{equation*}
  (H^2-G)dA
\end{equation*}
is invariant under conformal changes of the ambient metric.  If
$M$ is compact of genus $g$, then the Willmore energy
\begin{equation*}
  W = \int_M (H^2-G)dA+2\pi\chi(M)
\end{equation*}
is also conformally invariant. Note that if $\bar K=0$ and $M$ is
compact then the Willmore energy is the $L^2$--norm of the mean
curvature, because $2\pi \chi(M)=\int_M K dA=\int_M G dA$.

Compact Bryant surfaces with smooth ends in $\ball\subset \R^3$ and
complete minimal surface of finite total curvature with planar ends in
$\R^3$ extend, adding a finite number of ends, to compact surfaces in
the conformal 3--sphere $S^3=\R^3\cup \{\infty\}$. Their Willmore
energy is related to the intrinsic absolute total curvature with
respect to the corresponding space form geometry by
\begin{equation*}
  \int_{M\setminus\{ \text{ends}\}} |K|dA= \int_M (H^2-G)dA=W-2\pi \chi(M).
\end{equation*}
The first equality follows from the Gauss equation, because $H^2=-\bar K$ for
both Bryant surfaces and minimal surfaces in Euclidean space.

\begin{Pro}\label{Pro:compact}
  An immersed compact surface in $\R^3$ whose non--empty intersection
  with $\ball$ is a Bryant surface is a compact Bryant surface with
  smooth ends. In particular, except for finitely many ends the
  immersion is contained in $\ball$.
\end{Pro}

\begin{proof}
  Let $f\colon M \to \R^3$ be a conformal immersion of a compact
  Riemann surface $M$ and let $M_0= f^{-1}(\ball)$ such that
  $f_{|M_0}\colon M_0 \to \ball$ parametrizes a complete Bryant
  surface. The pull back to $M_0$ of the hyperbolic metric on $\ball$
  has non--positive curvature $K$, because $K=G+\bar K= -
  (H^2-G)\leq0$, and  the total absolute curvature of $f_{|M_0}$ is
  bounded by $\int_M (H^2-G)dA$.  So $f_{|M_0}$ induces a
  complete metric with non--positive curvature and bounded total
  curvature on $M_0$, which implies that $M\setminus M_0$ is a finite
  number of points, cf.~\cite[Chapter~III, Proposition~16]{La80}.
\end{proof}


\end{document}